%
\input amstex 
\documentstyle{amsppt}
\input bull-ppt
\keyedby{bull372e/mhm}

\define\TRACE{\operatorname{Trace}}



\define\TG{\widetilde{G}}
\define\HG{\widehat G}

\define\GL{\operatorname{GL}}
\define\CL{\operatorname{Cl}}
\define\SL{\operatorname{SL}}
\define\SO{\operatorname{SO}}
\define\AD{\operatorname{Ad}}
\define\Irr{\operatorname{Irr}}


\define\pie{\widetilde{\pi}}

\define\vp{\varphi}

\topmatter
\cvol{28}
\cvolyear{1993}
\cmonth{April}
\cyear{1993}
\cvolno{2}
\cpgs{215-252}
\title 
An external approach to unitary representations
\endtitle
\shorttitle{An external approach to unitary representations}
\author Marko Tadi\'c \endauthor
\address Department of Mathematics, University of Zagreb,
Bijeni\v{c}ka 30, 41000 Zagreb,  Croatia \endaddress
\cu Sonderforschungsbereich 170,
Geometrie und Analysis,
Mathematisches Institut Bunsenstr. 3-5, D-3400 
G\"{o}ttingen, Germany \endcu
\ml tadic\@cfgauss.uni-math.gwdg.de \endml

\endtopmatter

\document

\ah Introduction
\endah

The  principal ideas of harmonic analysis on a locally  
compact group $G$ which
is  not necessarily  compact or commutative were
developed in the 1940s and
 early 1950s. In this theory, the role of the classical 
fundamental
harmonics is played by the irreducible unitary 
representations of Ê$G$.
The set of all equivalence classes of such representations 
is denoted by 
$\HG$ and is called the dual object of Ê$G$Ê or the 
unitary dual of $G$. 

Since the 1940s, an intensive study of the foundations of 
harmonic analysis on
complex and real reductive groups has been in progress 
(for a definition of
reductive groups, the reader may consult the appendix at 
the end of \S 2).  The
motivation for this development came from mathematical 
physics, differential
equations, differential geometry, number theory, etc. 
Through the 1960s,
progress in the direction of the Plancherel formula for 
real reductive groups
was great, due mainly to Harish-Chandra's monumental work, 
while at the same
time, the unitary duals of only a few groups had been 
parametrized. 

With 
Mautner's work  [Ma], a study of harmonic analysis on 
reductive groups over
 other locally compact nondiscrete fields was started. We 
shall first
describe such fields. In the sequel, a locally compact 
nondiscrete field will
be called a local field. 

If we  have a nondiscrete absolute value on the field 
Ê$\Bbb Q$ Êof  rational
numbers, then it is equivalent either to the standard 
absolute value (and the
completion is the field Ê$\Bbb R$Ê of real numbers) or to 
a $p$-adic
absolute value for some prime number Ê$p$. For $r \in \Bbb 
Q^\times$ Êwrite
Ê$r=p^\alpha a/b$  where Ê$\alpha,a$, and $b$ are integers 
and neither $a$
nor $b$ are divisible by $p$. Then the $p$-adic absolute 
value of $r$ is 
$$ |r|_p=p^{-\alpha}.
$$
A completion of $\Bbb Q$
with respect to the $p$-adic absolute   value is denoted 
by $\Bbb Q_p$. It
is called a field of $p$-adic numbers. Each 
finite-dimensional extension
$F$ of $\Bbb Q_p$ has a natural topology of a vector space 
over $\Bbb
Q_pÊ$. With this topology, F becomes a local field. The 
topology of F can 
also be introduced with an absolute value which is denoted 
by  $|\; \;|_F$
Ê(in \S 5 we shall fix a natural absolute value). The fields
of real and complex numbers, together with the finite 
extensions of $p$-adic
numbers, exhaust all local fields of characteristic zero 
up to isomorphisms
[We].  

Let $\Bbb F$ be a finite field.
Denote by $\Bbb F((X))$ the field of formal  power series 
over $\Bbb F$. 
Elements of this field are series of the form  
$f=\sum_{n=k}^\infty a_n X^n, 
a_n \in
 \Bbb F$, for some integer $k$. Fix $qÊ>Ê1Ê$. Very  often 
$q$ is taken to be
the cardinal number of the finite field $\Bbb F$. One 
defines an  absolute
value of $f$ by the formula
$$
|f|_{{\Bbb F} ((X))}=q^{-\text{min}\{n;a_n \ne 0\}}
$$
when $f \ne 0$. In this way $\Bbb F((X))$ becomes a local  
field. Fields of formal  power series
over finite fields exhaust all local  fields of positive 
characteristic up
to isomorphisms [We].

The fields $\Bbb R$ and $\Bbb C$ are
called archimedean fields.
For any $x,y
\in \Bbb R^\times$ or $x,y \in \Bbb C^\times$, there is 
always a positive
integer $n$ such that  $|y| < |nx|$.   The  above property 
does not hold for any other
local field. This is the reason that local fields which 
are not isomorphic
to $\Bbb R$ or $\Bbb C$ are called nonarchimedean local 
fields. 

After Mautner in
the 1960s, a series of people started to consider 
reductive groups over
nonarchimedean local fields. Let us recall that $p$-adic 
fields were
introduced historically  to enable one to consider a 
single equation over a
$p$-adic field instead of an infinite series of 
congruences mod $p^k$. Arithmetical problems also provided 
motivation to consider representations of
reductive groups over such fields. The strongest 
motivation comes from the
LanglandsÕ program. A unifying element in this program is 
the representation theory of
reductive groups.  A nice introduction to the LanglandsÕ 
program is 
[Gb3].

Let $G$ be a reductive group over a local field.  
Harish-Chandra created a
strategy for obtaining the unitary dual $\HG$ through the 
nonunitary  dual
$\TG$, where $\TG$ is the set of all functional 
equivalence classes of
topologically completely irreducible continuous 
representations of $G$.
Functional equivalence means that the matrix coefficients 
of one representation
may be approximated by matrix coefficients of another on 
compact sets, and vice
versa. A complete definition of $\TG$ is in  \S 2. To 
obtain $\HG$, one needs to
classify $\TG$ (the problem of the nonunitary dual)  and 
to identify $\HG
\subseteq \TG$ (the unitarizability problem). In  [L2]  
Langlands showed how to
parametrize $\TG$ by irreducible representations with 
certain good  asymptotic
properties (tempered representations) of reductive 
subgroups, when the field
$F$ is $\Bbb R$. The tempered representations were 
classified for $F=\Bbb R$ by
Knapp and Zuckerman [KnZu] on the basis of 
Harish-Chandra's work, thus 
providing a complete picture of $\TG$. Despite the 
Langlands classification of
$\TG$ in the archimedean case, there were no big 
breakthroughs in  the
classification of unitary duals for quite a long time. 
Borel-Wallach and
Silberger proved that Langlands parametrization of 
$\widetilde{G}$ in terms of
tempered representations of reductive subgroups was valid 
for reductive groups
over all local fields [BlWh, Si1].

In this paper, we shall be concerned with the  
unitarizability problem for
reductive groups over local fields. One usually breaks the 
unitarizability
problem into two parts. The first part is constructing 
elements of
$\HG$, and the second is showing that the constructed 
representations
exhaust $\HG$ (completeness argument). The completeness 
argument is usually
realized by showing that the classes of  $\TG \backslash 
\HG$ are not
unitarizable. We may call such  an approach to the 
completeness argument
indirect. 

Suppose that the field is archimedean. Then one can 
linearize the
problems for  $\HG$ (and $\TG)$; one can ``differentiate'' 
the
representations  and come to the  infinitesimal theory 
where the main
object is the Lie algebra $\germ g$ of $G$. Also, for a 
maximal compact
subgroup $K$ of $G$, the   theory of compact Lie groups 
gives an explicit
description of $\widehat{K}$. Thus, one may try to 
understand $(\pi,H) \in \TG$
by studying the restriction of 
$\pi$
to $K$. These two points explain  why, in the  archimedean 
case, some
problems concerning representations, and  especially the 
unitarizability
problem, were often approached by studying the internal 
structure of
representations. Let us recall that the internal approach 
was very
successful in the compact Lie group case (restriction to a 
maximal torus).
In the nonarchimedean case, there is no possibility of 
such an internal
approach to the unitarizability problem. One of the 
reasons that there has
been much less study of the unitarizability problem is 
that the nonunitary
duals are not yet completely parametrized there. 

Despite the fact that unitary  duals of a very restricted 
number of groups have
been classified, it is interesting to note that in 1950 
Gelfand and  Naimark
published a book  [GfN2] in which they constructed what 
they assumed to be the
dual objects of the complex classical simple Lie groups. 
Their lists were very
simple, and the representations were also simple (although 
infinite
dimensional). Gelfand and Naimark were using functional 
analytic methods as
tools in their analysis. In 1967 Stein constructed, in a  
fairly simple manner,
representations in $\GL (2n,\Bbb C) \hat{\ }$  which were 
not contained in the
lists of Gelfand and Naimark [St]. For some other 
classical groups, it was even
easier to see the incompleteness of the lists from  [GfN2]. 

The
representations of Gelfand and Naimark  of $\GL (n,\Bbb 
C)$, complemented by
Stein, were not generally expected to exhaust the whole of 
$\GL (n,\Bbb
C) \hat{\ }$. 

The main aim
of this paper is to present  the ideas which lead first to 
the solution of
the unitarizability problem for $\GL (n)$ over 
nonarchimedean local fields
[Td3] and to the recognition that the same result holds 
over  archimedean
local fields [Td2], a result which was proved by Vogan 
[Vo3] using
an internal approach. Let us say that the  approach that 
we are going to
present may be characterized as external. At no point do 
we go into the
internal structure of representations. 

Let us present the answer. We fix a general  local field 
$F$. Let $D^u$ be the 
set of all functional equivalence classes of irreducible 
square integrable
modulo center representations of all $\GL (n,F), n \geq 1$ 
(for definition see
\S 3). Let $|\; \; |_F$ be the modulus of $F$ (see \S 5). 
For each representation  $\delta \in D^u$ of $\GL (n,F)$, 
and for
each $m \geq 1$, consider the representation of $\GL 
(mn,F)$ parabolically
induced (\S 1) by  
$$
|\det |_F^{(m-1)/2} \delta \otimes |\det |_F^{(m-3)/2} 
\delta 
\otimes \cdots  \otimes |\det |_F^{-(m-1)/2} \delta 
$$
 from a suitable
standard parabolic subgroup (i.e., from  one containing 
the upper triangular
matrices, \S 2). This representation has a unique  
irreducible
quotient which will be denoted by $u(\delta,m)$. For 
$0<\alpha<1/2$ let  $\pi
(u(\delta,m),\alpha)$ be the representation of $\GL 
(2mn,F)$ parabolically induced
by  $$
|\det |_F^\alpha u(\delta,m)  \otimes |\det |_F^{-\alpha}
u(\delta,m). 
$$
Denote by $B$ all
possible $u(\delta,m)$ and  $\pi (u(\delta,m),\alpha)$. 
Then the answer is

\proclaim{Theorem} {\rm (i)} Let 
$\tau_1, \ldots , \tau_n \in B$.
Then the representation 
$\pi$ parabolically induced by
$$
\tau_1 \otimes \cdots  \otimes \tau_n
$$
of suitable $\GL (m,F)$ is irreducible and unitary.

{\rm (ii)} Suppose that
$\rho$
is obtained from
$\sigma_1, \ldots , \sigma_n \in B$
in the same manner as
$\pi$
was obtained from
$\tau_1, \ldots , \tau_m$
in {\rm (i)}. Then
$\pi \cong \rho$
if and only if $n=m$ and the sequences
$(\tau_1, \ldots , \tau_n)$
and
$(\sigma_1, \ldots , \sigma_n)$
coincide after a renumeration.

{\rm (iii)} Each irreducible  unitary representation of 
$\GL (m,F)$, for
any $m$, can be obtained as in {\rm (i)}. 
\endproclaim

A new, and at the same time very old, point of
view that led to the papers  [Td1--Td3]
was that the unitarizability  problem has a reasonably 
simple  answer and
that the unitary representations appear in simple and 
natural ways. In the
completeness argument, instead of an indirect strategy, a 
direct argument was
used. In this way, a detailed study of nonunitary 
representations was
avoided. In all arguments, essentially only Hilbert space 
representations
were necessary.

Surprisingly, the statement  of the theorem, which was 
first discovered in
the nonarchimedean case in  [Td3], says that in the case 
$F=\Bbb C$ the
unitary dual of $\GL (n,\Bbb C)$ should consist of the 
representations of
Gelfand, Naimark, and Stein. Not only the statement of the 
nonarchimedean
case of the theorem, but also the methods of the proof in  
[Td3] made sense in
the archimedean case. Thus, after a complete proof had 
been written in the
nonarchimedean case, we wrote in  [Td2]
the proof of the archimedean case of the theorem.  We have 
used  there a
theorem of Kirillov from [Ki1], for which he never 
published the
complete proof (see \S 9).

There are now many solutions of the unitarizability 
problems,  especially for
particular reductive groups  in the archimedean case. In 
general, they are
based on ideas different  than the one that we present in 
this paper. Two of
them take a distinguished place (they solve completely the 
problem for a
series of groups having no bounds on their semisimple 
split ranks). The
first is Vogan classification in [Vo3] of the unitary 
duals of $\GL (n)$
over $\Bbb R, \Bbb C$, and $\Bbb H$. He proved a theorem 
equivalent to the
statement of our theorem for archimedean $F$. The other 
one is
Barbasch's classification  of the unitary duals of the 
complex classical
Lie groups in [Bb]. 

Since we consider both the archimedean and  nonarchimedean 
cases, it is
natural to recall  Harish-Chandra's Lefschetz principle: 
``Whatever is true
for real reductive groups is also true for $p$-adic 
groups'' [Ha2]. One
problem with the Lefschetz principle is that we  usually 
obtain a result
for archimedean $F$ by one kind of considerations and for 
nonarchimedean
$F$ by very different methods, and after that we compare 
the results. An interesting
problem is to explain the phenomenon of the Lefschetz 
principle, which is
certainly related to our depth of understanding of these two
theories. An important task of this paper is to present a 
unified point of
view on the theorem and its proof. We shall discuss both 
the theorem and
the proof, making no distinction concerning the nature of 
the field. This is
possible by the use of the external point of view. In our 
approach we
shall very often be close to the point of view of the 
representation theory
of general locally compact groups, and this will be a 
unifying point. 

Professor P. J. Sally suggested that I write a paper where 
the ideas would be
extracted  from the technical machinery as much as 
possible. I am thankful to
him for his advice and encouragement and for his generous 
help during the
writing of various drafts of this paper, which he has 
read. I am also thankful
for the hospitality and support of the University of 
Chicago during the summer
of 1986 when a sketch of this paper was prepared. I would 
like to express my
thanks to a number of mathematicians for helpful 
discussions or encouragement
during the writing of this paper and the preceding ones on 
the 
same topic. Let me mention just a few of them: M. Duflo, 
I. M. Gelfand, H.
Kraljevi\'c,
D. Mili\v ci\'c, F. Rodier, P. J. Sally, and M. F. 
Vigneras.  H.
Kraljevi\'c gave very helpful comments on the first draft 
of this paper.
Discussions with D. Mili\v ci\'c and his suggestions were 
of  great help in the
preparation of the paper. The referee gave a number of 
valuable and helpful
suggestions to the preceding draft of this paper. Among 
others, the work of J.
Bernstein and A. V. Zelevinsky greatly influenced the 
development of some ideas
presented in this paper. The final step in preparation for 
 publishing this
paper was done in G\"ottingen. I am thankful for  
hospitality and support to
the Sonderforschungsbereich 170 and
S. J. Patterson.

We hope that this paper will illustrate a certain internal 
symmetry in the 
external approach to the unitarizability problem. We hope 
that some of our
ideas will be helpful in dealing with the unitary duals of 
other nonarchimedean
classical groups (and also nonunitary duals). The papers  
[SlTd] and [Td8]
indicate that this hope is not without basis. 

Finally, we introduce some general notation which we shall 
use throughout the
paper. For  a topological space $X$,  $C(X)$ will denote 
the space of all
continuous (complex-valued) functions on $X$. The subspace 
of all compactly
supported continuous functions will be denoted by 
$C_{\roman c}(X)$. If we have
a measure $\mu$ on $X$, then $L^1(X, \mu)$ will denote the 
space of all classes
of $\mu$-integrable functions on $X$ and $L^2(X, \mu)$ 
will denote  the Hilbert
space of all classes of square integrable functions on $X$ 
with respect to
$\mu$. For a smooth manifold $X$ the space of smooth 
functions on $X$ will be
denoted by $C^\infty (X)$ and $C^\infty (X) \cap C_{\roman 
c}(X)$ will be
denoted by $C_{\roman c }^\infty (X)$. If $X$ is a totally 
disconnected locally
compact topological space, then $C_{\roman c}^\infty (X)$ 
will denote the space
of all compactly supported locally constant functions on 
$X$. The fields of
real and complex numbers are denoted by $\Bbb R$ and $\Bbb 
C$ respectively. The
ring of integers is denoted by $\Bbb Z$, nonnegative 
integers are denoted by
$\Bbb Z_+$, and positive ones are denoted by $\Bbb N$.


\ah 1. Concept of harmonic  analysis \\on general locally 
compact
groups
\endah

 In this section we shall   outline some of the ideas of 
harmonic analysis
on locally compact groups. 

Let $G$ be a locally compact group. We shall always 
suppose in this  paper that
the  groups are separable. A {\it representation\/} of $G$ 
is a pair $( \pi
,V)$  where $V$ is a complex vector space which is not 
zero dimensional and
$\pi$ is a homomorphism  of $G$ into  the group of all 
linear isomorphisms of
$V$. By a {\it continuous representation\/} of $G$ we 
shall mean a
representation $( \pi ,H)$ of $G$ where $H$ is a Hilbert 
space and the map
$(g,v) \mapsto \pi (g)v, G \times H \rightarrow H$ is 
continuous (we shall
always  assume that $H$ is separable). A closed subspace 
$H'$ of $H$ will be
called a  {\it subrepresentation\/} of a continuous  
representation $(\pi,H)$
of $G$, if $H'$  is invariant for all operators $\pi (g)$, 
$g \in G$. A
continuous representation $(\pi , H)$ of $G$ is called 
{\it irreducible\/} if
there does not exist a nontrivial subrepresentation of 
$(\pi ,H)$ (i.e.,
different from $\{0\}$ and $H)$.  A continuous 
representation $( \pi,H)$ of $G$
will be called a {\it unitary representation\/} if all the 
operators  $\pi(g),g
\in G$, are unitary. Two  unitary representations $(\pi_i, 
H_i), i=1,2$, of $G$
are called {\it unitarily  equivalent\/} if there exists a 
Hilbert space
isomorphism $\vp :H_1 \rightarrow H_2$ such that $ 
\pi_2(g)\vp= \vp  \pi_1(g)$ 
for all $g \in G$. The set of all unitarily equivalence 
classes of  irreducible
unitary representations of G will be denoted by $\HG$ and 
called the {\it 
unitary dual\/} of $G$. For a family $(\pi_i,H_i), i \in 
I$, of unitary
representations  of $G$, there is a  natural unitary 
representation of $G$ on
the direct sum of Hilbert spaces $\bigoplus _{i\in I}H_i$. 
This
representation will be denoted by $(\bigoplus _{i
\in I}\pi_i,
\bigoplus _{i\in I }H_i)$. It is called a {\it direct 
sum\/} of
representations $(\pi_i,H_i)$, $i \in I$.

A continuous representation (resp.\  unitary representation)
on a one-dimen\-sional space  is called a {\it 
character\/} (resp.\  {\it 
unitary character\/}) of $G$.

 The main problem of harmonic  analysis on the group $G$ 
is to understand
some interesting unitary representations of $G$ (such 
representations are
usually given on function spaces). One way to study this 
problem is to
break it into two parts (at least for type I groups which 
will be described
in the sequel of this section and which will be the only 
groups considered
in this paper). The first part is to
understand irreducible unitary representations, i.e., 
$\HG$, and the second
part is to understand other unitary representations in 
terms of irreducible
ones. This strategy is in the spirit of Fourier's 
classical idea of
fundamental harmonics. 

Let us explain what we mean by understanding general 
unitary  representations
in terms of irreducible ones. If $G$ is a compact group,   
then a fundamental
fact is that each unitary representation of $G$ can be 
decomposed into a direct
sum of irreducible unitary representations  [Di, Theorem 
15.1.3.].
Understanding a unitary representation $\pi$ in terms of 
irreducible unitary
representations means, in the compact case, to know how to 
decompose $\pi$ into
a direct sum of irreducible unitary representations. In 
the  noncompact case
each unitary representation decomposes into a direct 
integral of irreducible
unitary representations, and understanding here again 
means to know how to
decompose a given unitary representation into a direct 
integral of irreducible
representations. We are not going to define here the 
notion of direct integral
of representations because the definition is quite 
technical. The interested
reader may consult [Di, \S 8]. Let me only mention that 
the direct integral
generalizes the notion of direct sum and that direct 
integrals are determined
by measures on  $\HG$. To consider  measures on  $\HG$, 
one needs some
$\sigma$-ring of sets. This $\sigma$-ring arises in a 
standard way from the
natural topological structure on $\HG$. Now we shall 
define this topology. 

If $(\pi,H)$ is a continuous  representation of $G$ and 
$v,w \in H$, then the
function $g \mapsto (\pi(g)v, w)$ on $G$ is called a {\it 
matrix coefficient\/}
of $(\pi,H)$. We denote by $\Phi(\pi)$ the linear span of 
all matrix
coefficients of $(\pi,H)$.  The  closure  operator on 
subsets of $\HG$ is
defined as follows. Let $\pi \in \HG$ and $X \subseteq 
\HG$. Then   $\pi \in
\CL (X)$ if and only if each element of  $\Phi(\pi)$ can 
be approximated
uniformly on each compact subset of $G$ by elements  from $
\bigcup _{\sigma \in
X }\Phi(\sigma)$. This topology on $\HG$ will be called  
the {\it 
topology of the unitary dual\/} of $G$.

If $G$ is a commutative group,  then each irreducible 
unitary representation of
$G$ is given on a one-dimensional space (this follows 
easily from the spectral
theorem). Thus each  $\pi \in \widehat{G}$ is a function. 
Pointwise
multiplication of element of  $\HG$ defines a group 
structure on $\HG$.
Together with the above topology, $\HG$ is again a locally 
compact  commutative 
group which is  called the {\it  dual group\/} of $G$. The 
role of the topology
of the unitary dual is crucial in harmonic analysis on 
locally compact
commutative groups. This topology is a basis on which one 
builds the
fundamental  facts of the harmonic analysis on commutative 
groups. One of these
fundamental facts is that $G$ and $  (G\hat{\ }) \hat{\ }$ 
are canonically
isomorphic (Pontryagin duality). 

In the case of noncommutative groups, the role of this 
topology on $\HG$
is less crucial than in the  commutative case, but it is 
still an important
and natural object to consider. For example, $G$ is a {\it 
type I group\/} if
and only if $\HG$ is a $T_0$-space, i.e., for any two 
different points in
$\HG$, at least one of them has a neighborhood which does 
not contain the
other one [Di, 9.5.2. and Theorem 9.1]. It is important to 
notice that,
in general, $\HG$ is not topologically homogeneous and 
there exist 
significant connections between properties of irreducible 
representations
and their position in $\HG$ with respect to the topology.

Before we proceed further, we shall  say a few words about 
 some measures
which are natural to consider on locally compact groups. 
For a locally
compact group $G$, there exists a positive  measure 
invariant for right
translations. Such a measure will be called a {\it right 
Haar measure\/} on $G$
and it will be denoted by $\mu_G$. Thus,
$$
\int_G  f(gx) \, d\mu_G(g) \; =\; \int_G  f(g) \, d\mu_G(g)
$$
 for any $f \in C_{\roman c}(G)$ and $x\in
G$. Any two right Haar  measures are proportional.  The 
behavior of a
right Haar measure for left translations is described by 
the {\it 
modular function\/}. There exists a function $\Delta_G$ on 
$G$ such that
$$
\int_G  f(xg) \, d\mu_G(g) \;  =\; \Delta_G (x)^{-1} \; 
\int_G  f(g) \,
d\mu_G(g) 
$$
 for any $f
\in
C_{\roman c}(G)$ and $x
\in
G$. The group $G$ is called {\it unimodular\/} if
$\Delta_G \equiv 1$, i.e., if
$\mu_G$ is also invariant for left  translations. For  
more information
about Haar measures and for proofs of the above facts one 
may consult
[Bu2].

Suppose that $G$ is unimodular. 
The space $C_{\roman c}(G)$ becomes an algebra for the 
convolution
which is defined by the formula
$$
(f_1 * f_2)(x)=\int_Gf_1(xg^{-1})f_2(g)\, d\mu_G(g),
$$
$f_1, f_2 \in C_{\roman c}(G)$. In a natural way one can 
extend
the convolution to $L^1(G, \mu_G)$. Then $L^1(G, \mu_G)$
becomes a Banach algebra.
 For $f \in C_{\roman c}(G)$ and a
continuous representation $\pi$ of $G$
set
$$
\pi (f) \; =\;  \int_G  f(g) \pi(g) \, d\mu_G(g).
$$
Now
$\pi$ 
becomes a representation  of the convolution algebra  
$C_{\roman c}(G)$, and
this representation is called the {\it integrated form\/} 
of the representation
$\pi$ of $G$. If $\pi$ is unitary, then the last formula 
also defines a
representation of the algebra $L^1(G, \mu_G)$. Moreover, 
it is a
$*$-representation if we define 
$f^*(g)=\overline{f(g^{-1})}$. 

We have already mentioned that the  basic problem  of 
harmonic analysis on
$G$ is to classify $\HG$
and then to decompose interesting representations in terms 
of
$\HG$. Among the interesting  representations, there is 
one  that should be the
first to be understood, namely, the regular representation 
on $L^2(G,
\mu_G)$. By the abstract
Plancherel theorem, for a unimodular group $G$  there 
exists a unique positive
measure $\nu$ on $\HG$ such that
$$
\int_G  |f(g)|^2 \, d\mu_G(g)  \; =\; \int_{\HG}  
\text{Trace} ( \pi (f)
\pi (f)^*
 )\, d\nu(\pi) 
$$
for all $f \in L^1(G, \mu_G) \cap L^2(G, \mu_G)$ [Di,  
Theorem 18.8.2].
The measure $\nu$
is called the {\it Plancherel  measure\/} of $G$ (it  
determines an explicit
decomposition of $L^2(G, \mu_G)$ into a direct integral of 
elements of
$\HG)$.

While the basic ideas of  harmonic analysis  on general 
locally compact
groups were laid down in the 1940s, one of the first 
breakthroughs in
classifying unitary duals was the work of Kirillov for 
nilpotent Lie
groups at the beginning of the 1960s (see [Ki2, \S\S 13 
and 15]).

Let us first recall that a {\it Lie group\/}  is a group 
supplied with  a 
structure of a (real) analytic manifold such that the 
group operations
(i.e.,  multiplication and inversion) are analytic 
mappings. We can define
the {\it Lie algebra\/} $\germ g$ {\it of a Lie group\/} 
$G$ as the tangent
space of $G$ at the identity,  supplied with a bracket 
operation $[\;,
\;]$
which can be defined in the following way. If $X,Y \in 
\germ g$
are tangent vectors to curves $x(t), y(t)$ for $t=0$ 
respectively, then 
$[X,Y]$
is the tangent vector to the curve
$$
t \mapsto x(\tau) y(\tau) x(\tau)^{-1} y(\tau)^{-1}
$$
where
$\tau= \text{sgn} \, (t) |t|^{1/2}$,
at $t=0$ [Ki2, 6.3.]. An element $g
\in
G$ acts on $G$ by inner automorphism. The  differential  
of this action is
denoted by $\text{Ad}(g)$. In this way $G$ acts on $\germ 
g$, 
and this action is called the {\it adjoint action\/} of 
$G$ on $\germ g$. 

Let $G$ be a
connected simply connected nilpotent Lie group, and let 
$\germ g ^*$
be the space of linear forms on the Lie algebra
$\germ g$ of $G$. There is a natural action of $G$ on
$\germ g^*$.
It is called the {\it coadjoint action\/} of $G$. Using 
the theory of  induced
representations (which we  shall discuss a bit later in 
this section),
Kirillov established a canonical one-to-one 
mapping from the set of all coadjoint orbits onto the 
unitary dual of $G$
$$
G \backslash \germ g^* \rightarrow \HG
$$
which gives a simple description of
the unitary dual. With a natural    topology on the 
left-hand side, this is
a homeomorphism. Kirillov theory also gives the characters 
of irreducible
unitary representations and the Plancherel formula for 
nilpotent Lie groups. 

In the second part of this section, we
define some notions that we shall need in the sequel.

If $G$ is a compact group,
then each $\pi \in \widehat{G}$ is given on a  
finite-dimensional space [Di,
15.1.4.]. As in the theory of finite group 
representations, the function 
$$
\Theta_\pi:g \mapsto \TRACE \pi(g), 
$$
which is called the {\it character of the 
representation\/} $\pi$, completely
determines the class of $\pi$ in $\widehat{G}$. We have a 
right Haar measure on
$G$.  Thus  the character function determines a 
distribution on a compact Lie
group $G$. It is easy to see that  this distribution is 
$$
f \mapsto \TRACE \pi(f), f \in C^\infty _{\roman c}(G),
$$
 which will be denoted by
$\Theta_\pi$
again. If $\pi$ is an infinite-dimensional representation,
obviously the trace as a function is not well defined.
Nevertheless, it may happen that the above distribution
is well defined.
 For an arbitrary Lie group $G$,  one can  take the above 
distribution
for the definition of the character of $\pi \in \HG$, if
$\pi(f)$ is trace-class for 
$f \in C^\infty _{\roman c}(G)$.
If
$\pi(f)$ has a trace for
$f \in C^\infty _{\roman c}(G)$, then
$\pi
(f)$ must be a compact operator for any $f
\in
L^1(G,\mu_G)$. If all 
$\pi \in \HG$
have characters in the above sense, then 
$\pi(f)$ is a compact operator for any 
$\pi \in \HG$ 
and $f \in L^1(G,\mu_G)$.

A locally compact group is called a {\it CCR-group\/} if 
$\pi(f)$ is a compact operator for any
$\pi\in\HG$ and $f \in L^1(G,\mu_G)$.
 All CCR-groups are of type I [Di, Proposition 4.3.4. and 
Theorem 9.1]. 
A  great number of very important groups are CCR-groups. 
The most important
classes of CCR-groups are commutative groups, compact 
groups, nilpotent Lie
groups, and reductive groups over local fields (reductive 
groups will be
defined in the appendix at the end of the following 
section). In particular,
classical groups over local fields are CCR-groups. One can 
characterize
CCR-groups in terms of the topology of the unitary duals. 
A group $G$ is a
CCR-group if and only if $\HG$ is a $T_1$-space, i.e., if 
the points are closed
subsets of $\HG$ [Di, 9.5.3.].

For some important classes of groups  one  can show that 
they are
CCR-groups by showing that they have so-called ``large'' 
compact subgroups.
Now we shall explain the last notion. 

Let $G$ be a locally compact unimodular group, and let $K$ 
be a
compact subgroup. For $\delta \in \widehat K$ and a 
continuous representation 
$(\pi, H_\pi )$
of $G$, let $H_\pi(\delta)$
be the subspace of $H_\pi$
spanned by all subrepresentations of
$\pi|K$ which are isomorphic to
$\delta$
(here
$\pi | K$ denotes the  representation of $K$ obtained by 
restriction from
$G)$.  If 
$$
\text {dim}_{\Bbb C} H_\pi(\delta) < \infty
$$
for all
$\delta \in \widehat K$, then $\pi$ is called a {\it 
representation of\/} 
$G$ {\it with finite\/}
$K$-{\it multiplicities\/}. 

One  calls $K$ a {\it large compact subgroup\/} of $G$
if for each $\delta \in \widehat K$
the function
$$
\pi \mapsto \text{dim}_{\Bbb C}  \, \, H_\pi (\delta)
$$
is a bounded function on
$\HG$. If $G$ has a large compact subgroup, then $G$ is a 
CCR-group [Di,
Theorem 15.5.2]. Some very important classes of  groups 
have large compact
subgroups, for example, reductive groups over local fields. 

For a continuous
representation $(\pi,H)$ of $G$, a vector $v
\in
H$ is called $K$-{\it finite\/} if the span of all
$
\pi(k)v, k \in K,
$ is finite dimensional.

In the rest of this section  we shall discuss some  parts 
of the theory of
the induced representations for locally compact groups. 
The notion of
induced representations for locally compact groups 
generalizes the
well-known notion of induced representations for finite 
groups that was
introduced and studied by Frobenius and Schur. Induction 
is one
of the simplest and most important procedures for 
obtaining new
representations of locally compact groups. 

The most important case of induction for the purpose of 
this paper is {\it
para\-bolic induction\/}. To define this notion, it is  
enough to consider a
closed subgroup $P$ of a unimodular group $G$ and assume 
that there exists
a compact subgroup $K$ of $G$ such that $PK=G$. Let $( 
\sigma ,M)$ be a
continuous representation of $P$. The space of all 
(classes of) measurable
functions $f:G \rightarrow H$ which satisfy 
$$
f(pg)=\Delta_P(p)^{1/2} \sigma(p) f(g), \qquad p \in P, \ 
g \in G,
$$
 and 
$$
||f||^2=\int_K ||f(k)||^2 \; d\mu_K(k) \; < \; \infty
$$
 will be denoted by $\text{Ind}_P^G(\sigma)$. It is a
Hilbert space. The group $G$ acts by right  translations on
$\text{Ind}_P^G(\sigma)$. This action we denote by $R$. 
Thus 
$$
(R_gf)(x)=f(xg).
$$
With this action, $\text{Ind}_P^G(\sigma)$ is  a 
continuous representation
of $G$. It is unitary if $\sigma$
is unitary.

In our considerations $G$ will be a reductive group over  
a  local field, 
while $P $ will be a parabolic subgroup of $G$ (these 
terms will be defined in
the following section). One will take a Levi decomposition 
$P=MN$ of $P$ and a
continuous representation $\sigma$ of $M$. Since $P/N 
\cong M$, we shall
consider $\sigma$ as a representation of $P$. Then we 
shall say that
$\text{Ind}_P^G(\sigma)$ is a {\it parabolically induced 
representation\/} of
$G$ by $\sigma$ .

We shall also talk at some points in this
paper about induced representations which are of  more 
general type. Let
$G$ be a locally compact group which does not need to be 
unimodular, and
let $C$ be a closed subgroup of $G$. Suppose that a 
unitary representation $(
\sigma ,H)$ of $C$ is given (we may assume that  $\sigma$ 
is a continuous
representation only). We denote by 
$\text{Ind}_C^G(\sigma)$ the space of
all measurable functions $f:G \rightarrow H$ which satisfy  
$$
f(cg)=[\Delta_C(c) \Delta_G(g)^{-1}]^{1/2} \sigma(c) f(g), 
\qquad g \in G, \ c \in C.
$$
One square integrability condition  is required  also. 
This condition is
more technical than in the case of parabolic induction 
[Ki2, 13.2].
Again $G$ acts by right translations on 
$\text{Ind}_C^G(\sigma)$. This  
is a unitary representation of $G$. We say that 
$\text{Ind}_C^G(\sigma)$  is
{\it unitarily induced\/} by $\sigma$.

Mackey obtained a simple criterion for testing  if a given 
unitary
representation of $G$ is unitarily induced from $C$ 
(Imprimitivity Theorem,
[Ki2]). There is an important  specialization of this  
theory. If $G$
contains a nontrivial normal abelian subgroup $N$, then 
Mackey theory
implies a description of  $\HG$ by irreducible unitary 
representations of
smaller groups $(G$ and $N$ need to satisfy certain 
general topological
conditions). Let $\chi \in \widehat N$ and denote by 
$G_\chi$
the stabilizer of
$\chi$ in $G$ $(G$ acts on 
$\widehat N$
because $G$ acts on $N$ by automorphisms). Note that
$\widehat N$ consists of unitary characters. Take an 
irreducible unitary
representation $\sigma$ of $G_\chi$
such that
$\sigma | N$ is a multiple of
$\chi$. Then $\text{Ind}_{G_\chi}^G(\sigma)$ is an  
irreducible unitary
representation of $G$.  One obtains all irreducible
unitary representations of $G$ by the above construction. 
Two such
representations constructed  from $\chi_1, \sigma_1$ and
$\chi_2, \sigma_2$
are equivalent if and only if 
$\chi_1, \sigma_1$
and
$\chi_2, \sigma_2$
are conjugate. This specialization   of  Mackey theory is 
usually called
{\it small Mackey theory\/}. In the special case when $G$ 
is a semidirect
product of a closed normal abelian subgroup $N$ and a 
closed subgroup $M$,
the small Mackey theory describes $\HG$ more simply. Here 
$\HG$ is
parametrized by unitary characters of $N$ and irreducible 
unitary
representations of their stabilizers in $M$, divided by a 
natural
equivalence.

For our purpose, the case of semidirect  products is  the 
most interesting one.
It is important to observe that one obtains here 
automatically the
irreducibility of some induced representations. For a more 
detailed
exposition of small Mackey theory, one may consult [Ki2, 
13.3].


\ah  2. The nonunitary dual as a tool for the unitary dual
\endah

In the study of  the representation theory of the general 
linear groups, and
more generally, of the classical groups, the terminology 
of the 
theory of
reductive groups is very useful and natural. We shall very 
briefly
recall some of the terminology  of this theory in the 
first part of the
appendix at the end of this section. The reader may also 
skip over the general
definitions and follow only the case of $\GL (n)$ where we 
shall not need these
general definitions.

We begin  this
section with a few  definitions. 
We shall first define parabolic subgroups in $\GL (n,F)$. 
Take an ordered
partition  
$$
\alpha=(n_1,n_2,\ldots,n_k)
$$ 
of
$n$. Consider block-matrices 
$$
A =   
\left [
\matrix
A_{11}&.&.&.&.&.&A_{1k}\\
.&.&.&.&.&.&.\\
.&.&.&.&.&.&.\\
.&.&.&.&.&.&.\\
.&.&.&.&.&.&.\\
A_{k1}&.&.&.&.&.&A_{kk}
\endmatrix
\right]
$$
 where $A_{ij}$ is an $n_i$ by $n_j$ matrix. Denote 
$$
P_\alpha = \{ A \in \GL (n,F);A_{ij}=0 \; \; \text{for} 
\;\; i>j \},
$$
$$
M_\alpha = \{ A \in \GL (n,F);A_{ij}=0 \;\; \text{for} 
\;\; i \ne j \}.
$$
 Let
$N_\alpha$ be the set of all $A \in P_\alpha$  such that 
all $A_{ii}$ are
identity matrices. Now $P_\alpha$ is called a standard 
parabolic subgroup of
$\GL (n,F)$, $M_\alpha$ is called a Levi factor of 
$P_\alpha$, and $N_\alpha$ is
called the unipotent radical of $P_\alpha$. The subgroup 
$P_{(1,1, \dots, 1)}$
of all upper triangular matrices in $\GL (n,F)$ is called 
the standard minimal
parabolic subgroup.

Take any $g \in \GL (n,F)$ and any  ordered
partition  
$
\alpha=(n_1,n_2,\ldots,n_k)
$
of
$n$. Then $gP_\alpha g^{-1}$ is called a parabolic
subgroup of $\GL (n,F)$. Set $P'=gP_\alpha g^{-1}$, 
$M'=gM_\alpha g^{-1}$, and
$N'=gN_\alpha g^{-1}$. Then $P'=M'N'$ is called a Levi 
decomposition of $P'$.
The group $M'$ is called a Levi factor of $P'$, and $N'$ 
is called the unipotent
radical of $P'$. Similarly, a minimal parabolic subgroup 
is defined to be any
conjugate of the standard minimal parabolic subgroup.

We shall now introduce the nonunitary dual. We denote by 
$G$ a reductive group
over a local field $F$. There exists a maximal compact 
subgroup $K$ of $G$ such
that  $P_{\text{min}}K=G$  for some minimal parabolic  
subgroup
$P_{\text{min}}$ of $G$. We fix such a maximal compact 
subgroup. The Iwasawa
decomposition $ P_{\text{min}}K=G $ holds for any maximal 
compact subgroup $K$
of $G$ if $F$ is an  archimedean field. If it is not, this 
may not be true for
all maximal compact subgroups. The group $K$ is a large 
compact subgroup of $G$
in the sense of the previous section. In the case of $\GL 
(n,F)$ and $F=\Bbb R$
(resp.\ $F=\Bbb C)$, one may take for $K$ the group of 
orthogonal matrices
(resp.\ the group of unitary matrices). If $F$ is 
nonarchimedean, one may take
$K$ in $\GL (n, F)$ to be $\GL (n,\Cal O_F)$, where $\Cal 
O_F$ is the ring of
integers in $F$, that is, $\Cal O_F =\{ x \in F; |x |_ F 
\leq 1 \}$. In a
general linear group over any local field, all maximal 
compact subgroups are
conjugate.

We shall always assume in the sequel that continuous 
representations of $G$ that we consider have finite 
$K$-multiplicities.

For a continuous representation  $\pi$ of $G$, we have 
denoted by $\Phi(\pi)$
the linear span of all matrix coefficients of $\pi$. This 
is a subspace of
$C(G)$. Denote by $\CL  \Phi( \pi )$ the closure of $\Phi 
( \pi )$ with respect
to the open-compact \ topology on $C(G)$. Then two 
continuous
representations $\pi_1$ and $\pi_2$ of $G$ are said to be 
{\it functional
equivalent\/} if    
$$ \CL\Phi( \pi_1 )= \CL\Phi( \pi_2 ). $$ 
We denote by
$\TG$ the set of all functional equivalence classes of 
continuous  irreducible
representations of $G$ which have finite 
$K$-multiplicities. These
representations are precisely the topologically completely 
irreducible
representations of $G$ (for a definition of the last 
notion, one may consult
[Wr, 4.2.2.]). The set $\TG$ is called the {\it nonunitary 
dual\/} of $G$. We
could use Naimark equivalence to define $\TG$ instead of 
the functional
equivalence. We would get the same object.  Two continuous 
representations $(
\pi_1,H_1)$ and $( \pi_2,H_2)$ are called {\it Naimark 
equivalent\/} if there
exist dense subspaces $V_1 \subseteq H_1$ and $V_2 
\subseteq H_2$ which are
invariant for integrated  forms  and a closed one-to-one 
linear operator  $\vp$
from  $V_1$ onto $V_2$ such that   
$$ (\pi_2(f)\vp)(x)= (\vp  \pi_1(f))(x) $$
for any $f \in C_{\roman c}(G),x \in V_1.$

For an irreducible continuous representation
$\pi$
of $G$ which has finite $K$-multiplicities, the operator
$\pi(f)$ is of trace-class if
$f \in C_{\roman c}^\infty(G)$. The linear form
$$
f \mapsto \text{Trace}(\pi(f))
$$
is denoted by
$\Theta_\pi$ and called the {\it 
character of the representation\/}  $\pi$. Two  irreducible
continuous representations  with finite $K$-multiplicities 
are functionally
equivalent if and only if they have the same characters. 
Furthermore,
characters of representations from $\TG$ are linearly 
independent.

The natural mapping $\HG \rightarrow \TG $ is one-to-one 
$(G$ is a CCR-group).
Therefore, we shall identify $\HG$ with a subset of $\TG$. 
A class $\pi \in \TG
$ will be called {\it unitarizable\/} or {\it unitary\/} 
if $\pi \in \HG
\subseteq \TG$. One supplies $\TG$ with a topology in the 
same way as we have
supplied $\HG$ with the topology of the unitary dual.

The reader may consult the second part  of the appendix at 
the end of this
section for standard realizations of the set 
$\widetilde{G}$. Those realizations
depend on whether the field is archimedean or not.

The idea of Harish-Chandra was to break the problem of 
describing $\HG$ into
two parts: the  problem of the nonunitary dual and the 
unitarizability problem.
The problem of determining the nonunitary dual appeared to 
be much more
manageable than the unitarizability problem. 

Besides the above general strategy, there  could be other 
strategies for
getting $\HG$.  It would be interesting to obtain a 
classification of the
unitary duals dealing with nonunitary representations as 
little as possible  
or perhaps not at all. A strategy of such classification 
for $\GL (n)$  over
archimedean fields dealing with only unitary 
representations can be based on a
paper  from 1962 by Kirillov [Ki1]. We shall return to  
[Ki1] in \S 9. Now we
shall outline the strategy.

Let $P_n$ be the subgroup of $\GL (n,F)$  of all matrices 
with  the bottom  
row equal to $(0,0,\ldots ,0,1)$. An interesting property 
of $P_n$ follows
from the small Mackey theory: $\widehat{P}_n$ is in a 
bijection with
$$
\GL (n-1, F) \hat{\ }  \cup \GL (n-2, F)\hat{\ } \cup 
\cdots \cup \GL (2, F)
\hat{\ }  \cup \GL (1, F) \hat{\ } \cup \GL (0, F) \hat{\ }
$$
$(\GL (0,F)$ denotes the trivial  group). Gelfand and 
Naimark showed  already
in  1950 that certain irreducible  unitary representations 
of $\SL (n,\Bbb C)$,
which they expected to exhaust $\SL (n,\Bbb C) \hat{\ }$, 
remain irreducible as
representations  of $P'_n$ where $P'_n$ denotes the 
subgroup of $\SL (n,\Bbb
C)$ of all matrices with the bottom row of the form $(0, 
\ldots ,0,x), x
\in \Bbb C^\times$. They obtained this result from the 
explicit formulas
for those representations. Kirillov's idea was to first 
find a general
proof of the irreducibility of $\pi|P_n$ for $\pi \in
\GL (n,F)
\hat{\ }$ when $F=\Bbb R$ or $\Bbb C$. Then  one has an 
inductive procedure
for classification. After classifying $\GL (m,F)\hat{\ }$ 
for $m \leq
n-1$, one also has a classification of 
$\widehat {P}_n$ by the above remark about
$\widehat{P}_n$. Thus, the second part of the   strategy 
is to find all 
possible extensions of representations from  
$\widehat{P}_n$ to unitary
representations of $\GL (n,F)$. This strategy was used by  
I. J. Vahutinskii
in his study of irreducible unitary representations of 
$\GL (3,\Bbb R)$.

At the end of  this section we shall say a  few words 
about the characteristics
of some approaches to the unitary duals of certain groups 
in two relatively
simple cases. 

We shall first consider the case of a compact Lie group 
$G$. We shall assume
that $G$ is connected. One starts  from an  irreducible 
unitary representation
$(\pi ,H)$ of $G$, and studying  the internal structure of 
$H$, one comes to
exact parameters which classify $\HG$ . Let us give a 
rough  idea of how this
approach goes. A closed connected commutative subgroup of 
$G$ is called a {\it
torus\/} in $G$. Each torus is isomorphic to some $\Bbb 
T^n$ where $\Bbb T$ is
the group of all complex numbers of norm one. We fix a 
maximal torus $T$ in
$G$. Then each element of $G$ is conjugate to an element 
of $T$, i.e., each
conjugacy class of $G$ intersects $T$. Suppose that  $\pi$ 
is a continuous
representation of $G$ on a finite-dimensional space $H$. 
One can choose an
inner product on $H$ invariant for the action of $G$.
Thus the representation $\pi| T$ decomposes into a
direct sum of some unitary characters 
$\chi_1,\ldots,\chi_m$. These unitary
characters are called the {\it weights\/} of  $\pi$ with 
respect to $T$. Take
now $( \pi ,H) \in \HG$. We have already mentioned that 
the  character
$\Theta_\pi$ of the representation $\pi$ determines the 
class of $\pi$. Since
the character $\Theta_\pi$, as a function on $G$, is 
obviously constant on
conjugacy classes, $\pi$ is already determined by $\pi | 
T$, i.e.,  $\pi$ is
determined by its weights. Further analysis of  the 
structure of the
representation $\pi$ on $H$ requires the study of the 
representation  of the 
Lie algebra (the formula (L.A.) in the appendix at the end 
of this section,
defines the action of the Lie algebra). It gives that 
there is a particular
weight among all weights of $\pi$ which already 
characterizes $\pi$ (the
highest weight). In this way one obtains a  
parametrization of $\HG$ by a
certain subset of characters of $T$  (the dominant 
weights). For  an exposition
of this nice theory of Weyl and Cartan, one may consult, 
among many nice
expositions, the seventh paragraph of  [Bu1].

We shall present now a simple strategy  for solving the  
unitarizability problem for $G=\SL (2,\Bbb R)$. One may 
take for $K$ the
group $\SO (2)$ of all two-by-two orthogonal matrices of 
determinant one.
Note that $\SO (2) \cong T$. The unitary dual of $K$ is 
given by the
characters  $$
\sigma_n:
\left[
\matrix
\text{cos}(\vp) & -\text{sin}(\vp)\\
\text{sin}(\vp) & \text{cos}(\vp)
\endmatrix
\right]
\mapsto
e^{in\vp},
$$
when $n$ runs over $\Bbb Z$. If $(
\pi
,H)
\in \TG$, then it is not hard to show that  multiplicities 
of
the representation $\pi | K$ are one, i.e., $H(\sigma_n)$ 
are either zero- or
one-dimensional subspaces. This can be seen in a similar 
way as one shows that
the space of $K$-invariant vectors in an irreducible 
representation is one
dimensional if $(G,K)$ is a Gelfand pair [GfGrPi, Chapter 
III, \S 3, no. 4].
Thus,   there is a basis  $\{ v_n;n \in
S
\}$
of the Hilbert space $H$, parametrized by a subset $S$ of 
$\Bbb Z$.  Suppose
that $\pi$
is a unitary representation with a $G$-invariant  inner 
product
$(\;\;Ê,\;\;Ê)$ on $H$.  The formula (L.A.) from the 
appendix
defines the representation of the Lie algebra $\germ g$ of 
$G$ on
$K$-finite vectors. Differentiating the relation $( \pi 
(g)v, \pi
(g)w)=(v,w),g \in G$, along one-parameter subgroups,  one 
gets that the
representation on $K$-finite vectors satisfies
$$
(\pi(X)v,w)=-(v,\pi(X)w)
$$
 for all $X \in \germ g$. Put
$||
v_n||
=c_n$. Since $v_n,n
\in
S$, are orthogonal,  the inner product is  completely  
determined by
numbers $c_n, n \in
S$. One can solve the unitarizability problem in  the 
following way. Take
$( \pi
,H)
\in \TG$
and check if there exist positive numbers $c_n, n
\in
S$, such that the inner product
$$
\left(\sum \lambda_i v_i, \sum \mu_i v_i\right)=\sum c_i^2 
\lambda_i \bar{\mu_i}
$$
satisfies $(
\pi (X)v,w)+(v,
\pi
(X)w)=0$ for all $K$-finite vectors $v$ and $w$ in $H$ and 
all $X \in \germ g$. All
$\pi$ 
for which there exist such numbers  form the unitary dual. 
 Clearly, to
be able to solve the above problem, we need to know 
explicitly the
internal structure of irreducible representations of $G$.  

In the above two examples, one solves the problem
by study of the internal structure of representations.

\ah  Appendix \endah

\bh Algebraic groups \endbh
We shall recall very briefly in the  first part of this 
appendix some
definitions from the theory of algebraic groups. For 
precise definitions
containing all details, one should consult [Bl].

A {\it linear algebraic group\/} $\bold G$ is a Zariski  
closed  subgroup of
some general linear group with entries from an 
algebraically closed field. A
linear algebraic group is called {\it unipotent\/} if it 
is conjugate to a
subgroup of the upper triangular unipotent matrices. A 
linear algebraic group
$\bold G$  is called {\it reductive\/} if $\bold G$   does 
not contain a normal
Zariski closed unipotent subgroup of positive dimension. A 
linear algebraic
group $\bold G$   is called {\it semisimple\/}  if it is 
reductive and if it
has a finite center. If a noncommutative linear algebraic 
group  $\bold G$ does
not contain a proper normal Zariski closed subgroup of 
positive dimension, then
$\bold G$ is called a {\it simple algebraic group\/}. If 
there is a group
isomorphism of $\bold G$  onto some $\bold {\GL (1)}^n$ 
which is also an
isomorphism of algebraic varieties, then $\bold G$ is 
called a {\it torus\/}.

A Zariski closed subgroup $\bold P$ of $\bold G$ is called 
a {\it parabolic
subgroup\/} if the homogeneous space $\bold {G}/ \bold 
{P}$ is a projective
variety. 

By a {\it reductive group\/}, we shall mean the group 
$\bold  G (F)$ of all
$F$-rational points of a reductive group $\bold G$ which 
is defined as   an
algebraic variety over a local field $F$. By a {\it 
parabolic subgroup\/} of a
reductive group $\bold G (F)$, we shall  mean the group of 
all rational points
of a parabolic subgroup of $\bold G$, which is defined 
over $F$. A {\it minimal
parabolic subgroup\/} of $\bold G (F)$ is a parabolic 
subgroup which does not
contain any other parabolic subgroup. All minimal 
parabolic subgroups in $\bold
G (F)$ are conjugate. If we fix a minimal parabolic 
subgroup $P_{\text{min}}$
of $\bold G (F)$, then the parabolic subgroups containing 
$P_{\text{min}}$ are
called {\it standard parabolic subgroups\/}. Each 
parabolic subgroup in $\bold
G (F)$ is conjugate to a standard parabolic subgroup. Let 
$P$ be a parabolic
subgroup of $\bold G (F)$. Among Zariski connected normal 
unipotent subgroups
of $P$ there is the maximal one. It is called the {\it 
unipotent radical\/} of
$P$. Denote it by $N$. There exists a reductive subgroup 
$M$ of $P$ such that
$P=MN$ and $M \cap N = \{ 1 \}$. Note that $P$ is then a 
semidirect product of
$N$ and $M$. One says that $P=MN$ is a {\it  Levi 
decomposition\/} of $P$.
Also, one says that $M$ is a {\it Levi factor\/} of $P$.

Since a reductive  group $\bold G (F)$ is a closed  
subgroup of $\GL (n,F)$,
$\bold G (F)$ is in a natural way a locally compact group. 
If $F$ is an
archimedean field (i.e., $F=\Bbb R$ or $\Bbb C)$, then 
$\bold G (F)$ is a Lie
group in a natural way. If $F$ is a nonarchimedean field, 
then $\bold G (F)$ is
a totally disconnected group. 

In the sequel, a reductive  group $\bold G (F)$ will 
usually be denoted simply
by $G$.

The most important examples of reductive groups are the 
classical groups such
as general linear  groups $\GL (n,F)$, special linear 
groups, symplectic groups
and  orthogonal groups. The groups $\GL (n,F)$ form the 
simplest series of
reductive groups and one of the first series to be 
considered. 

Even if one is interested in harmonic analysis on some 
particular class of
classical groups, it is usually necessary to study a 
broader class of groups
because some important constructions involve subgroups 
which may not belong to
the considered class. Such subgroups are reductive. This 
is the reason why it
is convenient to use the terminology of reductive groups 
even if we consider
some specific class of groups. 


\bh Realizations of $\widetilde{G}$ \endbh
The set
$\TG$ has the following realizations, depending  on 
whether the  field is
archimedean or not. We shall not use these realizations in 
the sequel, so
the reader can also  skip over these definitions. We note 
that the
following notions are very important in the theory. It is 
also interesting
to note how different  $\TG$ looks in the following 
archimedean and
nonarchimedean realizations.

 Suppose that $F$ is
an archimedean field. First, we shall give  a definition 
of a $(\germ g
,K)$-module $(\germ g$ is the Lie algebra of $G)$. A {\it 
representation\/}
$\pi$ {\it of a Lie algebra\/} $\germ g$
is a real-linear map from
$\germ g$
into the space of all linear operators on a complex vector
space $V$ such that
$$
\pi([X,Y])=\pi(X)\pi(Y)-\pi(Y)\pi(X)
$$
for any $X,Y \in \germ g$. If we have a Lie group $G$ and 
a continuous
representation $\pi$
of $G$ on a finite-dimensional space $V$, then the 
following limit exists
$$
\pi(X)v\;=\; \frac{d} {dt} [\pi(x(t)v)]_{t=0},\qquad v \in 
V,
\tag{L.A.}
$$
and it defines a representation of the Lie algebra
$\germ g$
of $G$ on $V$ (in the above formula $X \in \germ g$
and $X$ is the tangent vector to  the curve $x(t)$ at 
$t=0)$. We call 
this Lie algebra representation the {\it differential\/} 
of $\pi$. Suppose
now that $G$ is a reductive group over $F$ and $\germ k$
the Lie algebra of $K$. Let $(\pi,V)$ be a pair where $V$ 
is a  complex
vector space and $\pi = (\pi_{\germ g}, \pi_K)$
is again a pair consisting of a Lie algebra representation 
$\pi_{\germ g}$ of
$\germ g$
on $V$ and  of a  representation $\pi_K$ of $K$ 
on $V$, such that the following three conditions are 
satisfied.
\roster
\item"{(a)}" {For each $v
\in
V$ the vector space $W$ spanned by all 
$\pi_K
(k)v, k
\in
K$, is finite dimensional and the  representation of $K$ 
on $W$ is
continuous. 
\item"{(b)}" The differential of the representation 
$\pi_K$ of $K$
equals the restriction of the Lie algebra representation
$\pi_{\germ g}$
to
$\germ k$.
\item"{(c)}" For any $k\in K, X \in\germ g$, 
and $v\in V$
$$
\pi_{\germ g}(\AD(k)X)v\;=\;\pi_K(k) \pi_{\germ g}(X) 
\pi_K(k^{-1})v. 
$$}
\endroster
Then $(\pi ,V)$ is called a $(\germ g ,K)$-{\it module\/}. 
An {\it
irreducible\/} $( \germ g ,K)$-{\it module\/} is a $(\germ 
g,K)$-module which 
has no nontrivial  subspaces invariant both for actions of 
$K $ and $\germ g$.
Two $(\germ g,K)$-modules $( \pi',V')$ and $( \pi'',V'')$ 
are {\it
equivalent\/} if there is a one-to-one linear mapping 
$\vp$ from $V'$ onto
$V''$ such that  $\vp \pi'_K(k)=\pi''_K(k) \vp $ and  $\vp 
\pi'_{\germ
g}(X)=\pi''_{\germ g}(X) \vp $ for any $k \in K$ and $X 
\in \germ g$. Now $\TG$
is in a natural bijection with the set of all  equivalence 
classes of
irreducible $(\germ g ,K)$-modules. If $( \pi ,H)$ is an 
irreducible continuous
representation of $G$ (with finite $K$-multiplicities, 
which we always assume),
then   one takes for $V$ the space of all $K$-finite 
vectors in $H$. The
formula (L.A.) defines an action of $\germ g$ on $V$, and 
$V$ becomes an
irreducible $( \germ g ,K)$-module in this way.

Suppose now that $F$ is nonarchimedean. A representation 
$( \pi ,V)$ of $G$ is
called {\it smooth\/} if for each $v \in V$ there is an 
open subgroup $K_v$ of
$G$ such that $ \pi(k)v=v$ for any $k \in K_v$. Again we 
say that smooth
representation $( \pi   ,V)$ is {\it irreducible\/} if 
there is no nontrivial
vector subspace  invariant  for the action of $G$. Two 
smooth representations
$(\pi_1,V_1)$ and $( \pi_2,V_2)$ of $G$ are {\it 
equivalent\/} if there exists
a one-to-one  linear map $\vp$ from $V_1$ onto $V_2$ such 
that  $\pi_2(g) \vp
=\vp \pi_1(g)$ for any $g \in G$. As before, there  is a 
natural one-to-one
correspondence from $\TG$ onto the set of all  equivalence 
classes of
irreducible smooth representations of $G$. If $( \pi  ,H)$ 
is an irreducible
continuous representation of $G$, one takes again for $V$ 
the space of  all
$K$-finite vectors $v$ in $H$. Now the restriction of the 
action of $G$ on $H$
to $V$ defines an irreducible smooth representation of $G$ 
on $V$.


\ah  3. Some simple constructions of unitary representations
\endah

One would like to have rather  simple and natural 
constructions  of unitary
representations which produce the whole of $\HG$. For 
nilpotent Lie groups such
a systematic procedure consists  of unitary induction by 
one-dimensional
unitary representations. For the groups we consider, the 
situation is not so
simple, but it is not too bad either. For example, one 
obtains the whole of
$\SL (n,\Bbb C) \hat{\ }$  by parabolic induction with 
one-dimensional, in
general nonunitary, representations (see \S 9). 

We have introduced the topology  of  $\HG$ in \S 1. In the 
construction of new
unitary representations,  the hardest problem is to find 
new connected
components of  $\HG$. Since $\HG$ is not topologically 
homogeneous,  there may
exist  special connected components, those consisting of 
only one point. These
representations are usually called {\it isolated 
representations\/}. To avoid the
influence of the commutative harmonic analysis coming from
$G^{\text{ab}}=G/G^{\text{der}}$ where $G^{\text{der}}$  
denotes the derived
group of $G$, we shall define the notion of 
representations isolated modulo
center. Let $Z(G)$ be the center of $G$. For $(\pi ,H) \in 
\TG$ there exists a
character $\omega_\pi \in Z(G)\tilde{\ }$ such that  
$\pi(z)= \omega_\pi (z) \,
\text{id}_H$ for all $z    \in Z(G)$. The character 
$\omega_\pi$ is called the
{\it central character\/} of $\pi$. For a character 
$\omega  \in Z(G)\tilde{\
}$ set 
$$
\TG_\omega=\{\pi\in \TG ; \; \omega_\pi =\omega \} \quad 
 \text{and} \quad 
\HG_\omega=\TG_\omega \cap \HG.
$$ 
The representation $\pi \in \HG$ (resp.\ $\pi \in \TG)$ 
will be called {\it 
isolated modulo center\/} (resp.\   {\it isolated modulo 
center in the nonunitary
dual\/}) if $\{\pi\}$ is an open subset of  
$\HG_{\omega_\pi}$ (resp.\ an open
subset of $\TG_{\omega_\pi})$. In the sequel, by isolated  
representation  we
shall mean isolated modulo center. According to what we 
have said about the
topology of $\HG$, we may say roughly that matrix 
coefficients  of isolated
representations are not similar to other matrix 
coefficients of elements of
$\HG$. The following example indicates that. If $G$ is 
compact, then $\HG$ is
discrete, and matrix coefficients of different  
representations are
$L^2$-orthogonal. Kazhdan proved in  [Ka] that the trivial 
representation is
isolated when $G$ is 
a simple group of split rank $n \geq 2$. The {\it split 
rank\/}
is the highest $n$ such  that $G$  possesses a Zariski 
closed subgroup defined
over $F$ which is  isomorphic over $F$ to $\GL (1,F)^n$. 
As opposed to the
trivial representation, other isolated representations are 
usually not easily
constructible. In fact, isolated representations of $\HG$ 
or $\TG$ are very
distinguished representations in known examples. 
Certainly, each isolated
representation in the nonunitary dual, which is unitary, 
is also isolated in
the unitary dual.

The
first isolated representations that one usually meets in the
representation theory of reductive groups 
are square integrable. An irreducible
unitary representation $(\pi ,H)$ of $G$ is called {\it 
square integrable modulo
center\/} if for any $v,w \in H$, the function  $$
g \mapsto |(\pi(g)v,w)|
$$
 is a square integrable function on
$G/Z(G)$ with respect to a Haar measure.  We shall use the 
term {\it square
integrable\/} instead of {\it square integrable modulo 
center\/}. Actually, the
unitarity of some $( \pi ,H) \in\TG$ may be obtained from 
the above
square integrability condition (note that the unitarity of 
the central character
of $\pi$ must be assumed to be able to formulate the above 
square integrability
condition). The square integrable representations are 
crucial for Plancherel
measure and important for parametrizing the nonunitary 
dual. In known examples,
they are very  often isolated (in the unitary dual).

 The known examples show that the construction of
a connected component, or at least a big part  of it, 
reduces to the
construction of isolated representations of reductive 
subgroups of $G$
attached  to parabolic subgroups and some standard simple 
and well-known
constructions.  Now we shall recall these standard simple 
constructions.
The first and the oldest one is:

\ch {\rm (a)} Unitary parabolic induction\endch
Let $P=MN$ be a Levi decomposition of a parabolic  
subgroup $P$  of $G$. For a
continuous representation $\sigma$ of $M$, we have 
considered $\sigma$ also as 
a representation of $P$ using the projection  $ P=MN 
\rightarrow P/N \cong M. $
Then $\text{Ind}^G_P(\sigma)$ was called a  parabolically 
induced
representation of $G$. If $\sigma$ is a unitary 
representation, then this
process  will be called unitary parabolic  induction. In 
fact, we always take
$\sigma \in \widehat{M}$. Then $\text{Ind}^G_P(\sigma)$ is 
a unitary 
representation which is a direct  sum of finitely  many 
irreducible
representations. It is usually irreducible. In the 
construction (a), we shall
always assume that P is a proper subgroup of $G$. In 
general, if
$(Ê\;\;,\;\;Ê)$ is an $M$-invariant hermitian form on the 
representation space
$U$ of $\sigma$, then  
$$
(f_1,f_2) \mapsto \int_K (f_1(k),f_2(k)) \; d\mu_K (k)
$$
is a $G$-invariant hermitian form on 
$\text{Ind}^G_P(\sigma)$, and it is
positive definite  if the form on $U$ was positive 
definite. 

Unitary parabolic induction has been used since  the first 
classification of
the unitary duals of reductive groups [Bg], [GfN1]. 
Gelfand and Naimark started
to use systematically   unitary parabolic induction for 
the classical simple
complex groups, while  Harish-Chandra started a systematic 
study of this
induction.

The following construction was used also
in the first classifications  of unitary duals  of reductive
groups [Bg, GfN1].

\ch{\rm (b)}  Complementary series \endch
It happens that some   representations  induced by 
nonunitary ones become
unitary after a new inner product is introduced on the 
representation space.
The idea is the following. One realizes on the same space 
a ``continuous''
family $(\pi_\alpha,H_\alpha), \; \alpha \in X$, of 
irreducible induced
representations  which possess $G$-invariant nontrivial 
hermitian forms. Let
$X$ be connected. Suppose that some $\pi_\alpha$ is 
unitary. The fact that a
continuous  family of nondegenerate hermitian forms on a 
finite-dimensional
space parametrized by $X$, being positive definite at one 
point of $X$, must be
positive definite everywhere enables one to conclude that 
all constructed
representations are in $\widehat{G}$ (here one reduces 
arguments to 
finite-dimensional spaces by considering spaces $\bigoplus 
H(\delta )$, where $\delta$
runs over fixed finite subsets of $\widehat{K})$. 
Positivity at  one point is
obtained in general from (a). The delicate point is the 
construction of
a continuous family of $G$-invariant hermitian forms, and 
it is based on the
theory of intertwining operators for induced 
representations. 

For the above construction some authors use the term {\it 
deformation\/} [Vo4].
We have chosen rather a more traditional name.

Let us recall that a topological space  $X$ is 
quasi-compact   if each open
covering of $X$ contains a finite open subcovering. Note 
that in the above
definition of a quasi-compact topological space, the 
Hausdorff property is not
required (this is the difference between quasi-compactness 
and compactness). A
topological space is locally quasi-compact if each point 
has a fundamental set
of neighborhoods which are quasi-compact. The fundamental 
fact about the
topology of $\HG$ (actually, about the dual of any 
$C^*$-algebra) is local
quasi-compactness. This fact essentially, together with 
some understanding of
the topology of the unitary dual, implies that 
$\widehat{G}$ cannot be complete
without  

\ch {\rm (c)}  Irreducible subquotients of
ends of complementary series \endch
This fact was first observed and proved by  
Mili\v{c}i\'{c}. Before we give a
brief  argument why the representations in (c) must be 
included in $\HG$, we
shall give a simple but suggestive example. 

Let
$P $ be a minimal parabolic subgroup  in $G=\GL (2,F)$ 
(one may take for
$P$ the upper triangular matrices in  $G)$. We have 
denoted by $\Delta_P$
the modular character of $P$. Representations 
$$
I^\alpha=\text{Ind}^G_P(\Delta_P^\alpha), \qquad -1/2 < 
\alpha < 1/2
$$
 are irreducible. If $\alpha=0$, then $I^0$ is unitary
by (a) $(\Delta_P^0$ is the trivial representation,  so it 
is unitary).
The family $I^\alpha,\, -1/2 <\alpha<1/2$, is a 
``continuous'' family of
irreducible representations with nondegenerate 
$G$-invariant hermitian
forms. Thus, they belong to $\widehat{G}$ by (b).

 We shall pay attention now to the representation at the end
of these complementary series 
$I^{-1/2}=\text{Ind}^G_P(\Delta_P^{-1/2}).$ From the 
definition of
$I^{-1/2}$, it follows that the trivial representation of 
$G$ is a
subrepresentation of $I^{-1/2}$ (the trivial   
representation of some
group $G$ is a representation on a one-dimensional space 
$V$ where each
element acts as the identity on $V)$. This is the unique 
proper nontrivial
subrepresentation of $I^{-1/2}$.  Since $I^{-1/2}$ is 
infinite
dimensional, there is no inner product on $I^{-1/2}$  for 
which
$I^{-1/2}$ is a unitary representation (in a unitary 
representation for
each subrepresentation there is another subrepresentation 
on the
orthogonal complement). Nevertheless, the representation 
on the
quotient of $I^{-1/2}$ by the trivial representation is 
unitary
(actually, it is square integrable, and this implies that 
it is unitary).
So, though $I^{-1/2}$  is not unitary, each irreducible   
subquotient of
the representation at the end of the complementary series 
is unitary.
This is the case in general. 

The representation $I^{-1/2}$ from the above  
considerations is a representation
which is not irreducible, but also, it is not very far 
from being irreducible.
This is an example of a representation of finite length. 
Before we proceed
further with explanation of the construction (c), we shall 
recall the
definition of a representation of finite length. Suppose 
that we have a
continuous representation $(\pi,H)$ of a reductive group 
$G$, which has finite
$K$-multiplicities. Then we say that $\pi$ is of {\it 
finite length\/} if there
exist subrepresentations  
$$  \{0\}=H_0 \subseteq H_1 \subseteq \dots \subseteq 
H_n=H $$
of $H$, such that the quotient representations of $G$ on 
$H_i/H_{i-1}$ are
irreducible representations of $G$, for $i=1,2,\dots,n$. 
Parabolic induction
carries the continuous representations of $M$ of finite 
length to the
continuous representations of $G$ of finite length.

Now we shall give a brief argument for the unitarity of 
representations in (c).
We shall omit  technical details. Suppose that we have a  
complementary series
of  representations $\pi_\alpha,\;\; \alpha \in X$. We may 
consider the
following situation. Let $Y$ be a topological space with a 
countable basis of
open sets, and let $X$ be a dense subset of $Y$. Assume 
that to each $\alpha
\in Y$ is attached a nontrivial continuous representation 
$\pi_\alpha$ of $G$
of finite length such that the functions  
$$ \alpha \mapsto \Theta_{\pi_\alpha}(\vp), \qquad Y 
\rightarrow \Bbb C $$ 
are continuous, for all $\vp$ from the  space $C_{\roman 
c,\ast}(G)$ of all
continuous compactly supported functions which span a 
finite-dimensional space
after translations by elements of $K$ (left and right). 
Suppose also that the
$\pi_\alpha$ are irreducible unitary representations for 
all $\alpha \in X$.
Take any $\alpha \in Y$. Let $(\alpha_n)$ be a sequence in 
$X$ converging to
$\alpha$. Mili\v{c}i\'{c} proved that in general a 
sequence $( \pi_n)$ in $\HG$
has no accumulation points  if and only if lim$_n \; 
\Theta_{\pi_n}(\vp) = 0$
for all  $\vp$ [Mi, Corollary of Theorem 6]. Since  
$\Theta_{\pi_\alpha}
\not\equiv 0$, $( \pi_{\alpha_n})$ has a convergent 
subsequence. If $S$ is the
set of  all limits of subsequences of $( \pi_{\alpha_n})$ 
in $\HG$, then
Mili\v{c}i\'{c}'s description of the topology  of $\HG$ 
says that there exist
positive integers $n_\sigma, \; \sigma \in S$, such that  
$$ \text{lim}_n \;
\Theta_{\pi_{\alpha_n}}(\vp) =  \sum_{\sigma \in 
S}\;\;n_\sigma
\Theta_\sigma(\vp) $$ 
for all $\vp$ in $C_{\roman c,\ast}(G)$ [Mi, Theorems 6 and
7].  Also, $S$ is a discrete and closed subset of 
$\widehat{G}$. Thus  $$
\Theta_{\pi_{\alpha}}(\vp) = \sum_{\sigma \in 
S}\;\;n_\sigma \Theta_\sigma(\vp)
$$ for all $\vp$ in $C_{c,\ast}(G)$. Since the character 
of  $\pi_\alpha$ is
the sum of characters of its  irreducible  subquotients, 
one can obtain easily
that $S$  is precisely the set of all irreducible 
subquotients of $\pi_\alpha$.
Since $S \subseteq \widehat{G}$, each irreducible 
subquotient of $\pi_\alpha$
is unitary. Thus (c) provides unitary classes. 

For a direct proof without use of the topology, one may 
consult 
[Td5]. The proof is based on the fact that a group of 
unitary operators on
finite-dimensional Hilbert space is finite dimensional. 
One uses in the
proof the fact that reductive groups have large compact 
subgroups.

While the constructions of (a) and (b) provide bigger   
continuous
families of unitary representations, (c) provides smaller 
families, but
they are often important in the construction of unitary 
representations.
Representations obtained by constructions (a), (b), or (c) 
are never 
isolated. We shall describe now a simple construction 
found by Speh
that may produce isolated representations. This 
construction is
particularly useful when it is combined with some other 
constructions,
for example, with (a), (b), and (c). Contrary to previous
constructions, here one gets unitarity of representations 
of smaller
groups from unitarity of representations of bigger groups. 
Before we
describe the construction, we need a notion of hermitian 
contragradient.

For a continuous representation $(\pi ,H)$ of $G$, 
$\bar{\pi}$ will denote the
{\it complex conjugate\/} of $\pi$.  It is the same 
representation, but the
Hilbert space is the complex conjugate of $H$. The {\it 
contragradient
representation\/} of $\pi$ will be denoted by $\pie$. It 
is the  representation
on the space of all continuous linear forms on $H$ with 
the action $[\pie
(g)f](v)=f(\pi(g^{-1})v)$. Set $\pi^+=\bar{\pie }$. Then 
$\pi^+$ will be called
{\it hermitian contragradient\/} of  $\pi$. A continuous 
irreducible
representation $\pi$ will be called {\it hermitian\/} if 
$\pi$ and $\pi^+$ are
in the same class in $\widetilde{G}$. It is easy to see  
that all unitary
representations are hermitian. In the classifications of  
$\TG$ it is usually
easy to check whether $\pi$ is hermitian or not.

Let $P=MN$ be a proper parabolic subgroup  of $G$. It  is 
very easy to
prove the following fact:

\ch{\rm (d)} Unitary parabolic reduction \endch
If we have hermitian
$\sigma \in \widetilde{M}$ such that 
$\text{Ind}^G_P(\sigma)$ is  
irreducible and that its class in $\TG$ is unitarizable, 
then the class
of $\sigma$ is unitarizable too.

We shall now give a rough argument  explaining why (d) 
provides  unitarizable
representations. Let $( \sigma ,H)$ be an irreducible 
nonunitarizable 
hermitian representation of $M$. Then there is a 
nondegenerate $M$-invariant
hermitian form $\psi$ on $H$. Now $H$ decomposes $H=H_+ 
\oplus H_-$ as a
representation of $K \cap P$, where $\psi$ is positive 
definite on $H_+$ and
negative definite on $H_-$.  Clearly, $H_+ \ne \{0\} $ and 
$H_- \ne \{0\} $.
Note  that $\text{Ind}^G_P(\sigma)$ and $\text{Ind}^K_{K 
\cap P}(\sigma| K \cap
P)$ are isomorphic as representations of $K$ and 
$$
\text{Ind}^K_{K \cap P}(\sigma| K
\cap P) \; \cong \; \text{Ind}^K_{K \cap P}(H_+) \oplus  
\text{Ind}^K_{K
\cap P}(H_-).
$$
Since unitary induction carries  unitary representations 
to unitary, we
see that the hermitian form induced by $\psi$
is indefinite. This form is also $G$-invariant.   Since a 
$G$-invariant
hermitian form on an irreducible representation is unique 
up to a scalar,
we see that if $\text{Ind}^G_P(\sigma)$ is irreducible, 
then it is not unitarizable.

Roughly speaking, the construction (d) enables one to   
construct sometimes
from a component of $\widehat{M}_1$ a component of 
$\widehat{M}_2$, where
$P_i=M_iN_i$ are two parabolic subgroups of $G$. 

There are also some simple constructions based on the 
geometry of a group or
groups. For example, if $\sigma_i \in \HG_i$, then 
$\sigma_1 \otimes \sigma_2
\in (G_1 \times G_2)\hat{\ }$ (and conversely). There are 
also irreducible
unitary   representations which appear already in the 
classification of the
nonunitary dual---the square integrable ones. 

It is interesting to ask which constructions must be  
added to (a)--(d) to
generate the whole unitary dual of the classical  groups, 
starting with
square integrable representations. We shall see that for 
the first class,
the case of $\GL (n,F)$, the constructions (a)--(d) are 
enough. 

\rem{Remarks} (1) Fell introduced in [Fe] a notion of 
nonunitary dual space 
for arbitrary locally compact group. It is a topological 
space consisting of
so-called linear system representations. We have studied 
$\TG$ as a topological
space  in [Td6] if $G$ is a reductive group over a 
nonarchimedean field $F$. 
It was shown that $\TG$ coincides with Fell's nonunitary 
dual. The set $\HG$ is
a closed subset of $\TG$. A representation $\pi \in \TG$ 
is isolated modulo
center if and only if there is a nontrivial matrix 
coefficient which is
compactly supported modulo center. Therefore, one may 
interpret Jacquet's
subrepresentation theorem [Cs, Theorem 5.1.2] in the 
following way. Each
element of $\TG$ can be obtained as  a subrepresentation of
$\text{Ind}^G_P(\sigma)$, with $\sigma$ isolated modulo 
center representation 
of $M$ for some parabolic subgroup $P=MN$ of $G$. Each  
isolated modulo center
representation of $G$ in $\TG$ is essentially unitary 
(i.e., it becomes unitary
after twisting by a suitable character of $G)$; actually 
it  is essentially
square integrable. Certainly, all these facts about the 
topology of the
nonunitary dual should hold over archimedean fields with 
the proofs along the
same lines. In the archimedean case, the representations 
with matrix
coefficients compactly supported modulo center can exist 
only if $G/Z(G)$ is
compact. Construction of representations of such groups is 
essentially  solved
by the case of compact Lie groups.

(2) Suppose that $F$ is nonarchimedean.  Let $P=MN$ be a  
parabolic subgroup in
$G$. Let $\sigma \in \widetilde{M}$ be isolated modulo 
center. Denote by $^0 G$
the set of all  $g \in G$ such that the absolute value of  
$\mu (g)$ is one for
any homomorphism $\mu :G \rightarrow F^\times$ which is 
also a morphism of
algebraic varieties defined over $F$. Then a character  
$\chi :G \rightarrow
\Bbb C^\times$ is called {\it unramified\/} if $\chi$ is 
trivial on $^0 G$,
i.e., if  $\chi(^0 G)=\{1\}$. Let $U(M)$ be the set of all 
 unramified
characters of $M$. Then $U(M)\sigma = \{ \chi \sigma;\; 
\chi \in U(M)\}$ is a
connected component of $\widetilde{M}$ containing  
$\sigma$,  and the set of
all irreducible  subquotients of $\text{Ind}^G_P(\tau), 
\tau \in U(M)\sigma$,
is a connected component of $\widetilde{G}$. All connected 
    components of
$\TG$ are obtained in this way [Td6]. So, for 
$\widetilde{G}$ the set of
connected components reduces  to  the set of isolated 
representations modulo
center in the nonunitary dual of the reductive parts of 
parabolic subgroups.
Note that a difficult problem in the nonarchimedean case 
is the construction of
representations isolated modulo center (i.e., of 
supercuspidal
representations).  
\endrem


\ah  4. Completeness argument 
\endah

In the last section, we have outlined  constructions 
(a)--(d) of unitary
representations of a reductive group $G$. It seems that 
those constructions
provide a  remarkable part of the unitary duals of the 
classical groups. This
leads to one of the most interesting questions about 
unitary representations of
reductive groups: How can one conclude that a set $X 
\subseteq \HG$ is a
significant piece of $\HG$, or even all of it? At the 
present time there is no
satisfactory strategy for answering such density 
questions.  Recall the simple
answer for finite groups: One needs to check if the sum of 
squares of degrees
of representations in $X$ is equal to the order of the 
group or not.

Suppose that a set $X \subseteq \HG$ is constructed and  
suppose also that one
expects that  it is the whole unitary dual. If one wants 
to prove that, then a
usual strategy has been to prove that in $\TG \backslash 
X$ all 
representations are nonunitarizable. One checks for each 
representation in $\TG
\backslash X$ that it cannot be unitarizable considering 
various properties of
that representation. The simplest properties that one can 
consider are: if the
representation is hermitian, if it has bounded matrix 
coefficients, etc. The
construction (d) can be used also for getting nonunitarity 
(from 
$\widetilde{M}$ to $\TG)$. The above strategy we shall 
call the indirect
strategy (of proving completeness of a given set of 
unitary representations). 

The indirect strategy becomes less satisfactory for groups 
of larger size. At
the same time, the indirect strategy is not completely 
satisfactory from the
point of view of harmonic analysis:  the stress is not on 
unitary
representations, which are of the principal interest, but 
on nonunitary ones. 
Actually, one needs a very detailed knowledge of the 
structure of
representations in $\TG \backslash \HG$, and the set $\TG 
\backslash \HG$ is
usually  much, much larger than the set $\HG$. The 
indirect strategy does not
develop directly the intuition about unitary 
representations.  In dealing with
$\TG$, it is very useful to algebracize the situation (the 
algebraic 
description of $\widetilde{G}$ was presented in the 
appendix at the end of \S
2).  

Later on we shall present a completeness argument for $\GL 
(n)$,  dealing
simultaneously with all $\GL (n)$ and having only one 
argument rather than
various cases. This will be an example of the direct 
strategy of proving
completeness.

In the following section we shall explain the sense in 
which the set of
representations of Gelfand, Naimark, and Stein is ``big''.


\ah  5. On the completeness argument: the example of
$\GL (n,\Bbb C)$
\endah

 We shall first introduce some general notation for the 
general linear group
and then pass to the complex case. 

In the first part of this section, $F$ denotes
any local field. 

For $x \in F^\times$, there exists a number $|x|_F>0$ such
that  
$$
|x|_F\; \int_Ff(xg)\;d_{\roman a}(g)=\int_Ff(g)\;d_{\roman 
a}(g)
$$
for all $f \in C_{\roman c}(G)$ $(d_{\roman a}(g)$ denotes 
an additive
invariant  measure on $F)$. Set $|0|_F=0$. Then $|\;\;|_F$ 
is called the
modulus of $F$. Note that $|\; \; |_{\Bbb R}$ is the usual 
absolute value on
$\Bbb R$,  $|\; \; |_{\Bbb C}$ is the square of the usual 
absolute value on
$\Bbb C$ (i.e.,  $|z|_{\Bbb C}=z\bar{z})$, while 
$|x|_{\Bbb Q_p}=|x|_p$ (see
the Introduction). For  $g \in \GL (n,F)$ set   
$$
\nu(g)=|\det (g)|_F.
$$
Clearly, $\nu:\GL (n,F) \rightarrow \Bbb R^\times$
is a character.

For $n_1,n_2
\in
\Bbb Z_+$, we have denoted by $P_{(n_1,n_2)}$ the 
parabolic subgroup  of
$\GL (n_1+n_2,F)$ consisting of the elements $g=(g_{ij})$ 
for which 
$g_{ij}=0$ when $i>n_1$ and $j \leq n_1$. Also we have 
denoted $M_{(n_1,n_2)}=
\{ (g_{ij}) \in
P_{(n_1,n_2)};g_{ij}=0\;\; \text{for}\;\; j>n_1 \;\; 
\text{and}\;\; i \leq
n_1\}$. Then 
$$
M_{(n_1,n_2)}
\cong
\GL (n_1,F)
\times
\GL (n_2,F)
$$
 is a Levi factor of $P_{(n_1,n_2)}$.

For two continuous representations 
$\tau_i$ of $\GL (n_i;F),\; i=1,2$, we  consider $\tau_1 
\otimes \tau_2$ as a
representation of $M_{(n_1,n_2)} \cong
\GL (n_1,F)
\times
\GL (n_2,F)$ in a natural way.
The mapping
$
mn \mapsto (\tau_1 \otimes \tau_2)(m),
$
where  $m \in M_{(n_1,n_2)}$ and $n \in N_{(n_1,n_2)}$, 
defines a 
representation of $P_{(n_1,n_2)}$. This representation of 
$P_{(n_1,n_2)}$ was
again denoted by $\tau_1 \otimes \tau_2$. Thus 
$$
\tau_1 \otimes \tau_2 :
\left[
\matrix
g_1&*\\
0&g_2
\endmatrix
\right]
\mapsto
\tau_1(g_1) \otimes \tau_2(g_2)
$$
for $g_1 \in \GL (n_1,F)$ and
$g_2 \in \GL (n_2,F)$.
Now the parabolically induced 
representation $$ \text{Ind}_{P_{(n_1,n_2)}}^{\GL (n_1+
n_2,F)}(\tau_1 \otimes
\tau_2) $$ will be denoted by
$
\tau_1 \times \tau_2$.
 It is a standard fact that
$(\tau_1 \times \tau_2) \times \tau_3$
is isomorphic to
$\tau_1 \times (\tau_2 \times \tau_3)$
(i.e., there exists a continuous intertwining which is  
invertible). 
This is a consequence of a general theorem on induction in 
stages 
$$
\text{Ind}_{H_2}^{H_3}(\text{Ind}_{H_1}^{H_2}(\sigma)) \cong
\text{Ind}_{H_1}^{H_3}(\sigma).
$$
 Therefore, it makes sense to write
$\tau_1 \times \tau_2 \times \tau_3$.

Before we explain an important property of the operation
$\times$, we need the notion of associate parabolic 
subgroups and  
associate  representations. Suppose that we have two 
parabolic subgroups
$P_1$ and $P_2$ of some reductive group $G$.  If we have 
Levi decompositions $P_1=M_1 N_1$, $P_2=M_2
N_2$, and $w \in G$ such that $M_2=wM_1w^{-1}$, then we 
say that $P_1$ and
$P_2$ are {\it associate parabolic subgroups\/}.  Suppose 
that $\sigma_1$
and $\sigma_2$ are finite length continuous 
representations of $M_1$ and
$M_2$ respectively, such that  $\sigma_2(m_2) \cong 
\sigma_1(w^{-1}m_2 w)$ for
all $m_2 \in
M_2$. Then we say that 
$\sigma_1$ and
$\sigma_2$ are {\it associate representations\/}.
For a continuous  representation $(\pi,H)$ of  $G$ of 
finite length, consider a
sequence of subrepresentations $$
\{0\}=H_0 \subseteq H_1 \subseteq \dots \subseteq H_k=H
$$
of $H$ where the quotient representations on $H_i/H_{i-1}$ 
are irreducible
representations of $G$, for $i=1,2, 
\dots,k$. Denote by $R(G)$ a free $\Bbb
Z$-module which has for a basis $\widetilde{G}$. We shall 
denote by  
$\text{J.H.}(\pi)$ the formal sum of all 
classes in $\widetilde{G}$ of the representations 
$H_i/H_{i-1},
\; i=1,2,\dots,k$. We shall consider $\text{J.H.}(\pi) \in 
R(G)$ and call it the
{\it Jordan-H\"older series of} $\pi$. The element 
$\text{J.H.}(\pi) \in R(G)$
does not depend on the filtration $H_i, \; i=0,1,2, \dots, 
k$ as above (one can
see that from the linear independence of characters of 
representations in
$\widetilde{G})$. Suppose that $P_1=M_1N_1$ and 
$P_2=M_2N_2$ are associate parabolic
subgroups. Let $\sigma_1$ and $\sigma_2$ be associate 
representations of $M_1$
and $M_2$ (as before, we consider $\sigma_1$ and 
$\sigma_2$ as
representations of $P_1$ and $P_2$ respectively). A 
standard fact  about 
parabolic induction from associate parabolic subgroups by 
associate
representations is that representations 
$\text{Ind}_{P_1}^G(\sigma_1)$ and
$\text{Ind}_{P_2}^G(\sigma_2)$ have the same characters, 
which implies that
these two  representations have the same Jordan-H\"older 
series. The formula for
the character of a parabolically induced representation 
from a minimal parabolic
subgroup, when $F=\Bbb R$, is computed in [Wr, Theorem 
5.5.3.1]. A similar
formula holds without assumption on parabolic subgroup. In 
the nonarchimedean
case we have a similar situation.

Suppose that $\tau_1$ and
$\tau_2$ are continuous representations of finite length 
of $\GL (n,F)$ 
and $\GL (m,F)$ respectively. The above fact about 
induction from
associated parabolic subgroups by associate 
representations implies that
$\tau_1 \times \tau_2$ and $\tau_2 \times \tau_1$
have the same Jordan-H\"older series.

Set  $\Irr ^u= \bigcup _{n \geq 0 }\GL (n,F)\hat{\ }$. To 
solve the
unitarizability problem for the $\GL (n,F)$-groups,  one 
needs to determine
$\Irr ^u$. 

In the rest of this section we shall  assume $F=\Bbb C$. 
Recall that
$|\;\;|_{\Bbb C}$ is the square of the standard absolute 
value that we
usually consider on $\Bbb C$.

Let
$\chi_0:\Bbb C^\times \rightarrow \Bbb C^\times$ be  the 
character $x \mapsto
x|x|_{\Bbb C}^{-1/2}$. Since each
$\pi \in
\GL (n,\Bbb C)
\hat{\ }$
has a central character and $\GL (n,\Bbb C)$ is a product 
of the   center
and of $ \SL (n,\Bbb C)$, the restriction of 
representations from
$\GL (n,\Bbb C)$ to $\SL (n,\Bbb C)$ gives a one-to-one 
mapping of  
$$
\GL (n,\Bbb C)\hat{\ }_{\chi_0} \cup \GL (n,\Bbb C)\hat{\ 
}_{\chi_0^2}  \cup
\cdots 
 \cup  \GL (n,\Bbb C)\hat{\ }_{\chi_0^n}
$$
 onto $\SL (n,\Bbb C)\hat{\ }$. Therefore, in order to 
understand 
$\SL (n,\Bbb C)\hat{\ }$ it is enough to understand $\GL 
(n,\Bbb C)\hat{\ }$
(and conversely). In the rest of this paper we shall deal 
only   with
$\GL $-groups and interpret the Gelfand, Naimark, and 
Stein representations in
terms of $\GL (n)$.

The first obvious irreducible unitary representations  of 
$\GL (n,\Bbb
C)\hat{\ }$ are one-dimensional representations 
$\chi \circ \det _n$
 where
$\chi$
is a unitary character of $\Bbb C^\times$ and where $\det 
_n$ denotes the
determinant  homomorphism of $\GL (n)$. Gelfand and 
Naimark obtained also the
following series of irreducible unitary representations 
$$
(\nu^{-\alpha} \chi) \times (\nu^\alpha \chi)    =
[\nu^{-\alpha}(\chi \circ \det _1)] \times 
[\nu^\alpha(\chi \circ
\det _1 )],\qquad \chi \in (\Bbb C^\times)\hat{\ },\  0 < 
\alpha <
\tfrac 12,
 $$
 the complementary
series representations for $\GL (2,\Bbb C)$. These 
complementary series 
start from representations 
$$
 \chi \times  \chi   =(\chi \circ
\det _1) \times (\chi \circ \det _1) ,\qquad \chi \in (\Bbb
C^\times)\hat{\ }.
 $$
 The following unitary representations of $\GL (n,\Bbb C)$ 
will  be
obtained by para\-bolic induction using the above 
representations.  Gelfand
and Naimark showed in [GfN2] that the unitary 
representations obtained by
parabolic  induction using the representations above are 
irreducible.
They assumed that in this way one gets all irreducible 
unitary
representations of $\GL (n,\Bbb C)$, up to unitary 
equivalence. 

We
can interpret the above remarks in the following way. Set 
$$
B_0=\{\chi \circ \det _n,\; [\nu^{-\alpha}(\chi \circ
\det _1)] \times [\nu^\alpha(\chi \circ \det _1 )];\;\chi  
\in
(\Bbb C^\times)\hat{\ }, \;  n \in \Bbb N, \; 0 < \alpha < 
\tfrac 12\}.
$$
Gelfand and Naimark
showed that for $\tau_1, \ldots ,\tau_k \in B_0$, the 
representations
$\tau_1 \times \cdots \times \tau_k$
are in $\Irr ^u$. Their assumption was that in this way 
one can get any 
representation from 
$\Irr ^u$.

Stein showed that the Gelfand and Naimark complementary 
series  
representations for $\GL (2,\Bbb C)$ are just the first of 
the
complementary series representations which exist for all 
$\GL (2n,\Bbb C)$
[St]. He showed that $$
[\nu^{-\alpha}(\chi \circ
\det _n)] \times [\nu^\alpha(\chi \circ \det _n )] \in 
\Irr ^u
$$
also for 
$n \geq 2, \; 0< \alpha < 1/2, \; \chi \in (\Bbb 
C^\times)\hat{\ }$, 
and he showed that these representations were not obtained 
 by Gelfand
and Naimark. These complementary series start from 
representations 
$$
(\chi \circ \det _n) \times (\chi \circ \det _n),\qquad n  
\geq
2,\ \chi \in (\Bbb C^\times)\hat{\ },
$$
 which were already
well known to Gelfand and Naimark. 

Now it is natural to complete the Gelfand and
Naimark  list by the Stein complementary series 
representations. 
Therefore, put  
$$
B=\{\chi \circ \det _n,\; [\nu^{-\alpha}(\chi \circ
\det _n)] \times [\nu^\alpha(\chi \circ \det _n )];\;\chi  
\in
(\Bbb C^\times)\hat{\ }, \;  n \in \Bbb N, \; 0 < \alpha < 
\tfrac 12\}
$$ 
(i.e., $B$ is just $B_0$ completed by the Stein 
complementary  series
representations). Now using arguments similar to those of  
Gelfand and
Naimark, one may conclude that for $\tau_1, \ldots ,\tau_k 
\in B$ the
representations $\tau_1 \times \cdots \times \tau_k$
are in $\Irr ^u$(see
[Sh]). Let us denote by 
\{G.N.S.\}
the set of all representations obtainable in this way.

We can explain now in what sense
\{G.N.S.\}
is big in $\Irr ^u$. It is easy to prove the following 
fact   (and we shall
prove it later): 

\thm\nofrills{(D)\usualspace} For any $\pi \in \Irr ^u$, 
there exist
$\tau_1,\;\tau_2 \in \{\roman{G.N.S.}\}$  such that  $\pi 
\times \tau_1$ and
$\tau_2$ have a composition factor in common.
\ethm

It is clear that (D) plays a role in the completeness 
argument.

Regarding unitary parabolic induction for $\GL (n,\Bbb 
C)$,  the simplest 
situation would be if it were always irreducible. This is  
what Gelfand
and Naimark expected to hold in 1950:

\thm\nofrills{(S1)\usualspace} 
Unitary parabolic induction for
$\GL (n,\Bbb C)$ is irreducible, i.e.,
$$
\tau_1,\; \tau_2 \in \Irr ^u \;\;\;\;\; 
\text{implies}\;\;\;\;\;\;  \tau_1
\times \tau_2 \in \Irr ^u.
$$
\ethm 

Let us suppose that (S1) holds. Because of (S1) and (D), 
for each $\pi \in \Irr
^u$ there exist  $\tau_1,\;\tau_2 \in$ \{G.N.S.\} such 
that  $\pi \times
\tau_1=\tau_2$. Thus, to obtain $\Irr ^u$, it is enough to 
 know how
representations from  \{G.N.S.\} can be parabolically 
induced.

We shall call $\pi \in \HG$ {\it primitive\/} if there is 
no proper parabolic
subgroup $P=MN$ and $\sigma \in \widehat{M}$ so that $\pi 
\cong
\text{Ind}_P^G(\sigma)$. Certainly, each $\pi \in \HG$ is 
unitarily equivalent
to some $\text{Ind}_P^G(\sigma)$ where  $\sigma \in 
\widehat{M}$ is primitive.
We have mentioned that if $P_1$ and $P_2$ are associate   
parabolic subgroups
and $\sigma_1,\;\sigma_2$ associate representations, then
$\text{Ind}_{P_1}^G(\sigma_1)$ and 
$\text{Ind}_{P_2}^G(\sigma_2)$ have the same
Jordan-H\"older series.  The simplest situation  would be 
if the converse were
true for $\sigma_1$ and  $\sigma_2$ primitive (because of 
the induction in
stages, it is necessary  to assume that $\sigma_1$ and  
$\sigma_2$ are
primitive). For $\GL (n,\Bbb C)$ this would mean  (having 
(S1) in mind): 

\thm\nofrills{(S2)\usualspace}
If \ both families $\tau_i \in \GL (n_i, \Bbb C)\hat{\ 
},\;\;
i=1,\dotsc,n,$ and $\sigma_j \in 
\GL (m_j,\Bbb C)
\hat{\ }, 
j=1,\dotsc,m,$ consist of primitive representations, and if
$$
\tau_1 \times \cdots \times \tau_n=\sigma_1 \times  \cdots 
\times
\sigma_m, $$
then $m=n$, and after a renumeration,  the sequences 
$(\tau_1, \ldots,
\tau_n)$ and  $(\sigma_1, \ldots$ 
$\ldots , \sigma_m)$ are
equal \RM(all $n_i,\;m_j$ are assumed to be $\geq
1)$.
\ethm

We shall assume that this holds.

A very plausible hypothesis is:

\thm\nofrills{(S3) \usualspace}
The Stein representations
$$
[\nu^{-\alpha}(\chi \circ
\det _n)] \times [\nu^\alpha(\chi  \circ \det _n )];\qquad 
\chi \in
(\Bbb C^\times)\hat{\ }, \   n \in \Bbb N, \  0 < \alpha < 
\frac 12 ,
$$
are primitive.
\ethm

Certainly, if we assume (S3), then all elements of B are 
primitive. 
Now it is obvious that (D), together with (S1), (S2), and 
(S3), implies
$$
\Irr ^u=\{\text{G.N.S.}\}.
$$

It remains to find a way to prove (S1), (S2), and (S3).    
Note that 
until now only (classes of) unitary representations were
necessary, and statements (S1), (S2), and (S3) are 
essentially analytic.
So, if one could prove (D), (S1), (S2), and (S3) dealing 
only with unitary
representations, one would have a classification of $\GL 
(n,\Bbb C)\hat{\
}$ completely in terms of unitary representations. 

As we shall see, (D) is easy to prove using noncomplicated
parts of the nonunitary theory. The following  strategy 
for proving (S2)
and (S3) simultaneously can be used. The set $\Irr ^u$ can 
be embedded in a
suitable ring which is factorial and where multiplication 
corresponds to
parabolic induction. Then one can prove that elements of B 
are prime or
close to being prime. This would imply (S2) and (S3). This 
strategy would
also include nonunitary theory (but again, not the 
complicated parts). 

Finally, we leave the discussion about (S1) for \S 9. Let 
us say that the first
ideas  for proving (S1) are due to Gelfand, Naimark, and 
Kirillov. 

In the following sections we shall elaborate in more 
detail the above strategy
and outline such a strategy for $\GL (n)$ over general 
local field $F$.


\ah  6. The
nonunitary dual of $\GL (n,F)$
\endah

In this section, we state some basic facts
about a parametrization of 
$\GL  (n,F)^\sim$. 

Besides $\Irr ^u$ which was
introduced in the previous section, we introduce 
$$
\Irr  =  \underset n \geq 0 \to \bigcup \GL  (n,F)^\sim.
$$
The set of all classes of square
integrable representations in $\Irr ^u$ of all 
$\GL (n,F), n \geq 1$,  will be
denoted by $D^u$. The set of all essentially square 
integrable representations
will be denoted by $D$. More precisely, 
$$
D=\{(\chi \circ \det )\ \delta;\,\chi  \in \GL  
(1,F)^\sim,\;\;\delta
\in D^u\}. $$
 For a set $X$, $M(X)$ will denote the set of all finite 
multisets
in $X$. These are all unordered $n$-tuples,  $n \in\Bbb Z_+
$. This is an
additive semigroup for the operation 
$$
(a_1,\dotsc,a_n)+
(b_1,\dotsc,b_m)=(a_1,\dotsc,a_n,b_1,\dotsc,b_m).
$$
 
By $W_F$
we shall denote the Weil group of $F$ if $F$  is 
archimedean and the
Weil-Deligne group in the nonarchimedean case. For the
purposes of this paper, it will not be essential to know 
exactly the definition of
$W_F$. We denote by $I$ the set of all classes of 
irreducible 
finite-dimensional representations of $W_F$. By the local 
Langlands conjecture for
$\GL (n)$ (which generalizes local class field theory), 
there should exist
a natural one-to-one mapping of $\Irr$  onto classes of 
semisimple
representations of $W_F$, i.e., onto $M(I)$
$$
\Irr  \rightarrow M(I).
$$
 Under such a mapping, $D$
should correspond to $I$. Thus, there should exist a  
parametrization of
$\Irr $ by $M(D)$. Let us write one such parametrization. 

Let $a=(\delta_1, \ldots , \delta_n) \in M(D)$. We can 
write 
$$
\delta_i=\nu^{e(\delta_i)}\delta_i^u,\qquad e(\delta_i)  
\in \Bbb
R,\ \delta_i^u \in D^u.
$$
 After a renumeration, we may assume  $e(\delta_1)\geq 
e(\delta_2) \geq
\cdots \geq e(\delta_n)$. The representation
$$ 
\lambda(a)=\delta_1 \times \cdots \times \delta_n
$$
has a unique irreducible quotient (possibly 
$\lambda(a)$ itself), whose class depends  only on $a$ 
[BlWh,
Jc1]. Its multiplicity in $\lambda(a)$ is 1. We shall 
denote this class
by $L(a)$. Now $$
a
\mapsto
L(a), \;\;\;\;\;\;\;\;M(D)
\rightarrow
\Irr 
$$
is a one-to-one mapping onto $\Irr $. It is a version of 
the  Langlands
parametrization of the nonunitary duals of $\GL 
(n)$-groups.  Certainly,
for the existence of such a parametrization, it is crucial 
that the
parabolic induction by square integrable representations 
is irreducible
for $\GL (n)$. 

The formula for
the hermitian contragradient in the Langlands 
classification becomes 
$$
L((\delta_1, \ldots , \delta_n))^+=L((\delta_1^+, \ldots , 
\delta_n^+)).
$$
 If we
write $
\delta_i=\nu^{e(\delta_i)}\delta_i^u (e(\delta_i) \in  \Bbb
R,\;\;\delta_i^u \in D^u),
$ then 
$$
\delta_i^+=(\nu^{e(\delta_i)}\delta_i^u)^+= 
\nu^{-e(\delta_i)}\delta_i^u. 
$$
Also, it is easy to show that
$$
\nu^\alpha L((\delta_1, \ldots ,  \delta_n))=L((\nu^\alpha 
\delta_1,
\ldots , \nu^\alpha \delta_n)),\qquad  \alpha \in \Bbb C.
$$

 Let $R_n$ be the free abelian group
with basis $\GL (n,F)^\sim$, i.e., $R_n=
R(\GL (n\!,F))$. For each finite-length continuous 
representation $\pi$
of $\GL (n\!,F)$\<, we have denoted by 
$\roman{J.H.}(\pi)$ its Jordan-H\"older  series
which is an
element of $R_n$. Set 
$$
R = \underset n \geq 0 \to \bigoplus \, R_n.
$$
 Now $R$ is a graded additive group which is free over 
$\Irr $. Let
$$
R_n \times R_m \rightarrow R_{n+m}
$$
be the $\Bbb Z$-bilinear mapping defined on the basis 
$\Irr $ by
$$
(\sigma, \tau) \mapsto \text{J.H.}(\sigma \times \tau).
$$
 This defines a
multiplication on $R$, which will be denoted again by 
$\times$. We  have
mentioned in the previous section that $\tau \times \sigma$
and
$\sigma \times \tau$ have the same Jordan-H\"older 
sequences. This 
just means the commutativity of  $R$. Since 
$(\tau_1 \times \tau_2) \times \tau_3 \cong \tau_1 \times 
(\tau_2 
\times \tau_3)$ (\S 5), $R$ is also associative. Certainly,
$\{L(a);\;\;a\in M(D)\}$
is a basis of $R$. It is a standard fact that
$\{\lambda(a);\;\;a\in M(D)\}$
is a basis of $R$ (i.e., standard characters form  a basis 
of  the group
of all virtual characters). This means nothing else than 
the following
fact which was first noticed by Zelevinsky in the 
nonarchimedean
case [Ze, Corollary 7.5]. 

\proclaim{\num{6.1.} Proposition} 
The ring $R$ is a $\Bbb Z$-polynomial  algebra over
all  essentially square integrable representations 
\RM(i.e., over $D)$.
\endproclaim

 In particular, $R
$ is a factorial ring. 

It is a natural question   to ask if it is
possible  to relate the operation of summing 
representations of $W_F$
with some operation on representations from $\Irr $ in the 
correspondence
$a \mapsto L(a)$; i.e., is there a relation between $L(a+
b)$ and
representations $L(a)$, $L(b)$? The answer is very 
nice:$L(a+b)$ is
always a subquotient of $L(a) \times L(b)$. 

One may find more information about the Langlands 
classification on  the level
of general reductive groups in  [BlWh]. For the Langlands 
philosophy one may
consult [Gb3].


\ah  7. Heuristic construction
\endah

In this section we shall try to see what would be   the 
part of the unitary
dual for the $\GL (n)$-groups over arbitrary local field 
$F$, generated by
classical constructions (a)--(d) of \S 3. Our principle 
will be
to expect a situation as simple as we could assume, 
bearing   our
evidence in mind. 

So, let
us start with $D^u$. The representation $\delta \times 
\delta$
is irreducible, so
$\delta \times \delta \in
\Irr ^u$ by the construction (a). Then we  have the  
complementary series
which starts from $\delta \times \delta$:
$$
[\nu^\alpha \delta] \times [\nu^{-\alpha}\delta],\qquad  0 
< \alpha <1/2.
$$
This is an example of construction (b).  We require    
$\alpha<1/2$ to
ensure that the induced representation is irreducible. At 
the end of
complementary  series $[\nu^\alpha \delta] \times 
[\nu^{-\alpha}\delta]$, 
there will be unitary irreducible subquotients.  To be 
able to identify
at least one, let us recall the relation mentioned at the 
end of \S 6; namely,
$$L(a+b) \text { {\it is a subquotient of\/} }L(a) \times 
L(b)
\text{ {\it for \/} } a,b \in M(D)\.$$
We shall assume that this holds in the rest of this section.
Therefore, we have that
$L((\nu^{1/2}\delta,
\nu^{-1/2}\delta))$ is unitary since it is at the end of   
the above
complementary series (construction (c)).

To proceed further, in order to be able to form a new
complementary series, let us suppose that for general 
linear groups
{\it  unitary
parabolic induction is irreducible}.
\noindent
 Then one has $L((\nu^{1/2}\delta,
\nu^{-1/2}\delta)) \times L((\nu^{1/2}\delta,
\nu^{-1/2}\delta)) \in \Irr ^u$, and further, one  has a 
new complementary
series 
$$ 
[\nu^\alpha L((\nu^{1/2}\delta,
\nu^{-1/2}\delta))] \times [\nu^{-\alpha}L((\nu^{1/2}\delta,
\nu^{-1/2}\delta))],\;\; 0 < \alpha <1/2.
$$
At the end of this complementary series
$$\align
&[\nu^{1/2} L((\nu^{1/2}\delta,
\nu^{-1/2}\delta))] \times [\nu^{-1/2}L((\nu^{1/2}\delta,
\nu^{-1/2}\delta))]\\
&\qquad =L((\nu \delta, \delta)) \times L((\delta, 
\nu^{-1}\delta)) 
\endalign
$$
is the unitary subquotient
$$
L((\nu \delta, \delta, \nu^{-1}\delta,  \delta))
$$
by the assumption that $L(a+b)$ is a subquotient of $L(a)
\times
L(b)$. It is natural to ask if  the above representation  is
primitive. But $
L((\nu \delta, \delta, \nu^{-1}\delta,  \delta))
$ is a subquotient of 
$$
L((\nu \delta, \delta, \nu^{-1}\delta)) \times L((  
\delta)).
$$
Note that $
L((\nu \delta, \delta, \nu^{-1}\delta)) \otimes L((  
\delta))
$ is hermitian. If we take the simplest possibility, that 
is, $
L((\nu \delta, \delta, \nu^{-1}\delta)) \times L((  
\delta))$ 
irreducible, then the construction (d) implies that
$$
L((\nu \delta, \delta, \nu^{-1}\delta)) \otimes L((  
\delta))
$$
is unitary. Clearly then, $
L((\nu \delta, \delta, \nu^{-1}\delta)) \in \Irr ^u
$.

At this point it is convenient to introduce some notation. 
Set
$$
a(\gamma,n)=(\nu^{(n-1)/2}\gamma,
\nu^{(n-3)/2}\gamma,\ldots,\nu^{-(n-1)/2}\gamma)
$$
 for
$\gamma \in
D$ and $n
\in \Bbb Z^+$ and
$$
u(\gamma,n)=L (a(\gamma,n)).
$$
We shall write
$$
\nu^\alpha(\delta_1,\ldots,\delta_n)=(\nu^\alpha \delta_1,
\ldots,\nu^\alpha \delta_n).
$$

 Now we proceed further. We already have 
$u(\delta,1),\;u(\delta,2),\;u(\delta,3) \in \Irr ^u$. One 
considers a
complementary series 
$$
[\nu^\alpha u(\delta,3)] \times [\nu^{-\alpha} 
u(\delta,3)],\;\;0<\alpha<1/2,
$$
which starts from $u(\delta,3)
\times
u(\delta,3)$. At the end we have the representation
$$
L(\nu^{1/2} a(\delta,3)) \times L(\nu^{-1/2} a(\delta,3)),
$$
and we can identify one irreducible subquotient which is 
$$
L(\nu^{1/2} a(\delta,3) + \nu^{-1/2} a(\delta,3)).
$$
It is a unitary subquotient by (c). Note that
$$
\nu^{1/2} a(\delta,3)+\nu^{-1/2} a(\delta,3)=a(\delta,4) +
a(\delta,2)
$$
which means that $L(a(\delta,4)+a(\delta,2))$ is unitary. 
Further,  this
representation is a subquotient of  
$$
u(\delta,4) \times u(\delta,2)=L(a(\delta,4)) \times 
L(a(\delta,2)).
$$
Suppose again that $u(\delta,4)
\times
u(\delta,2)$ is irreducible. Since it is unitary, 
$u(\delta,4)
\otimes
u(\delta,2)$ needs to be unitary by the construction (d). 
Therefore,
$u(\delta,4)$ will be unitary. 

Now it is easy to conclude that the assumption
$$
u(\delta,n) \times u(\delta,n-2 ) \in \Irr 
$$
leads to
$$
u(\delta,n)  \in \Irr ^u.
$$

In the case of $\GL (n,\Bbb C)$, $D^u$ is equal to $\GL 
(1,\Bbb C)\hat{\ }
=(\Bbb C^\times)\hat{\ }$. Then  
$u(\delta,n)=\delta
\circ
\det _n$. So, we have not obtained new unitary 
representations in 
this case.


\ah  8. Scheme of unitarity for $\GL (n)$
\endah

In the last section we have seen heuristically what some 
simple 
assumptions  suggest. Now we shall write down some of 
those assumptions
and their ``implications''. We first recall that 
$$
u(\delta,n)=L (a(\delta,n))=L((\nu^{(n-1)/2}\delta,
\nu^{(n-3)/2}\delta,\ldots,\nu^{-(n-1)/2}\delta))
$$
for $\delta \in D$ and $n \geq 1$.
It was also observed in \S 6 that $R$ is a factorial ring.

 We introduce the  following
statements: 

\roster
\item"{(U0)}"
$\tau,\;\sigma \in \Irr ^u \;\Rightarrow \; \tau \times 
\sigma \in \Irr $.

\item"{(U1)}"
$\delta \in D^u,\;\; n\in \Bbb N\;\Rightarrow \; 
u(\delta,n)  \in \Irr ^u.$

\item"{(U2)}"
$\delta \in D^u,\;\; n\in \Bbb N, \;\; 0 < \alpha < 1/2 
\;\Rightarrow \;
[\nu^{\alpha}u(\delta,n)] \times 
[\nu^{-\alpha}u(\delta,n)]   \in \Irr ^u.$

\item"{(U3)}"$\delta \in D,\;\; n\in \Bbb N\;\Rightarrow  \;
u(\delta,n)$ {is a prime element of} $R$.

\item"{(U4)}"$a,\;b\in M(D) 
\;\Rightarrow \; L(a+b)$
{is a composition factor of} $L(a)
\times
L(b)$.
\endroster

Here only (U3) was not assumed  or obtained in the last 
section.  It is
a strengthening of the assumption that the $u(\delta,n)$'s 
are primitive, 
which was present in the last section (otherwise, we would 
have tried to
construct  new unitary representations in that way). 

We have seen that (U0) and (U4) lead to
(U1) and (U2) (i.e., unitarity of the representations 
mentioned there).
 But it is interesting and surprising that with the 
addition of only one
assumption, namely, (U3), the preceding assumptions also 
easily imply
completeness for the unitary duals of the groups $\GL 
(n,F)$.

\proclaim{\num{8.1.} Proposition} Suppose that {\rm 
(U0)--(U4)} hold
true. Set 
$$
B=\{u(\delta,n),\;\;  [\nu^\alpha u(\delta,n)] \times 
[\nu^{-\alpha}
u(\delta,n)],\;\;\delta \in D^u,\;\;n \in \Bbb N,\;\; 0 < 
\alpha <\tfrac 12\}.
$$
 Then\/{\rm :}
\roster
\item"{(i)}" if $\tau_1, \ldots , \tau_k \in B$,
we have $
\tau_1 \times \cdots \times \tau_k \in \Irr ^u${\rm ;}
\item"{(ii)}" if 
$\pi \in \Irr ^u$, then there exist
$\sigma_1, \ldots , \sigma_m \in B$ such that 
$$
\pi=\sigma_1 \times \cdots \times \sigma_m;
$$
\item"{(iii)}" if
$\sigma_1, \ldots , \sigma_k,\;\;\tau_1, \ldots , \tau_m 
\in B$
and
$\sigma_1 \times \cdots \times \sigma_k=\tau_1 \times 
\cdots \times
\tau_m$, then  $k=m$ and the sequences 
$\sigma_1, \ldots , \sigma_k$ and $\tau_1, \ldots , \tau_m$
coincide after a renumeration.
\endroster
\endproclaim

From (U0), (U1), and (U2) one obtains  (i) directly. Also 
(U3) implies 
(iii). It remains to prove only (ii). 
First we shall prove 

\proclaim{\num{8.2.} Lemma} Suppose that
$\pi \in
\Irr $ is hermitian. If {\rm (U4)} holds, then there exist
$\sigma_1, \ldots , \sigma_n,\;\;\tau_1, \ldots , \tau_m 
\in B$
such that 
$\pi \times \sigma_1 \times \cdots \times \sigma_n$ and 
$\tau_1 \times 
\cdots \times \tau_m$
have a composition factor in common.
\endproclaim

\demo{Proof} Let $\delta \in
D^u,\;\; k
\in
(1/2)\, \Bbb Z_+$, and $ 0<\beta<1/2$. Then
$$
(\nu^{-k}\delta,\nu^{k}\delta) + 
a(\delta,2k-1)=a(\delta,2k+1),\qquad k>0,
$$
and
$$\align
&(\nu^{-k-\beta}\delta,\nu^{k+\beta}\delta) +
\nu^{\beta-1/2} a(\delta,2k) + \nu^{1/2-\beta}a(\delta,2k)\\
&\qquad =\nu^{\beta} a(\delta,2k+1) + 
\nu^{-\beta}a(\delta,2k+1).
\endalign
$$

Let
$\pi \in
\Irr $ be hermitian. Then
$\pi
=L((\gamma_1, \ldots , \gamma_s))$ for some $\gamma_i \in 
D$.  Now 
$$
L((\gamma_1, \ldots ,
\gamma_s))^+=L((\gamma_1^+, \ldots , \gamma_s^+))
$$ 
(see \S 6). This 
implies that we can write $\pi$ in the form
$$
\pi=L\left(\sum_{i=1}^{n_1}
(\nu^{-k_i}\delta_i,\nu^{k_i}\delta_i)+\sum_{i=n_1+1}^{n_2}
(\nu^{-k_i-\beta_i}\delta_i,\nu^{k_i+\beta_i}\delta_i)
+\sum_{i=n_2+1}^{n_3}(\delta_i)\right) 
$$
where $\delta_i \in
D^u,\;\; k_i \in
(1/2)\,\Bbb Z_+,\;\;  0<\beta_i<1/2$ and all 
$k_1,\ldots,k_{n_1}>0$. Here $n_1,
n_2, n_3 \in \Bbb Z_+$ and $n_2 \geq n_1, n_3 \geq n_2$ 
(i.e., not all three
sums need to show up in the above formula). The two 
relations at the beginning
of  the proof and (U4) imply that both 
$$ \pi \times \left[\prod_{i=1}^{n_1}
u(\delta_i,2k_i-1)\right] \times \left[\prod_{i=n_1+1}^{n_2}
(\nu^{\beta_i-1/2} u(\delta_i,2k_i)  \times
\nu^{1/2-\beta_i}u(\delta_i,2k_i))\right] 
$$
and
$$ \align  
\left[\prod_{i=1}^{n_1}
u(\delta_i,2k_i+1)\right] &\times
\left[\prod_{i=n_1+1}^{n_2}
(\nu^{\beta_i} u(\delta_i,2k_i+1) \times 
\nu^{-\beta_i}u(\delta_i,2k_i+1))\right] \\
&\times
\left[\prod_{i=n_2+1}^{n_3}u(\delta_i,1)\right]
\endalign
$$
have
$$ \align
L\left(\sum_{i=1}^{n_1}a \vphantom{\sum_{i=n_1+1}^{n_2}} 
\right.&
(\delta_i,2k_i+1)\\
&\left. +\sum_{i=n_1+1}^{n_2}
(\nu^{\beta_i} a(\delta_i,2k_i+1) + 
\nu^{-\beta_i}a(\delta_i,2k_i+1))
+\sum_{i=n_2+1}^{n_3}(\delta_i)\right) 
\endalign
$$
as a composition factor. This completes the proof of the 
lemma.\qed\enddemo

\demo{End of proof of Proposition \rm 8.1} Let
$\pi \in \Irr ^u$. Then, by the preceding lemma there exist
$\sigma_i$, $\tau_j \in
B$ so that
$\pi \times \sigma_1 \times \cdots \times \sigma_n$
and
$\tau_1 \times \cdots \times \tau_m$
have a composition factor in common. Since (U0), (U1),  
and (U2)  imply
that both sides are irreducible, we have  
$$
\pi \times \sigma_1 \times \cdots \times \sigma_n=\tau_1 
\times \cdots \times 
\tau_m. $$
Since $R$ is factorial and the $u(\delta,n)$'s are prime 
by (U3),
$\pi$
is a product of some $u(\delta,k)$'s, $\delta \in D$.  But 
the fact that
$\pi$ is hermitian implies that
$\pi$ is actually a product of elements of $B$. So, we 
have proved (ii).
\qed\enddemo

\rem{\num{8.3.} Remark} For $\GL (n)$ over a central 
simple division 
$F$-algebra, we expect that a scheme of this type should 
work too. In that case
a  slight modification in the definitions of the 
$u(\delta,n)$'s and  lengths
of complementary series is necessary (see [Td7]). 
\endrem


\ah  9. On proofs
\endah

We have seen that the fulfillment of (U0)--(U4) implies  a 
 complete
solution of the unitarizability problem for $\GL (n,F)$, 
as was described in
Proposition 8.1.

Let us remark that (U0)--(U4) were expected to hold for 
$F=\Bbb C$ 
 (except maybe (U3) because such questions were not 
considered).
Statements (U1) and (U2) were known, (U0) was expected 
even by Gelfand
and Naimark, while (U4) is easy to prove. A simple 
consequence of
(U0)--(U4), namely, the description of the unitary dual of 
$\GL (n,\Bbb C)$, 
was not generally expected to hold. Now we shall make a 
few remarks on
the history and proofs of (U0)--(U4). 

We start with the statement
(U4), which belongs to the theory of the nonunitary dual. 
This  fact was
proved by Zelevinsky in the case of nonarchimedean $F$ for 
his
classification of $\GL (n,F)$. His proof uses induction on
Gelfand-Kazhdan derivatives [Ze, Proposition 8.4]. Rodier 
noticed
in  [Ro] that Zelevinsky's proof implies (U4) for the 
Langlands
classification in the nonarchimedean case. We  proved (U4) 
in a  simple manner
for the archimedean case [Td2, Proposition 3.5. and 5.6]. 
Such a  proof is
outlined for nonarchimedean $F$ in Remark A.12(iii) of 
[Td3]. Sometimes one
can conclude the equality in (U4),  using the Zelevinsky's 
proof of
Proposition 8.5 in [Ze].

\proclaim{\num{9.1.} Proposition} Let 
$a_i=(\delta^i_1,\dotsc,\delta^i_{n_i})
\in
M(D),\;\; i=1,2$. If $\delta^1_k
\times
\delta^2_m
\in
\Irr $ for all $1
\leq
k
\leq
n_1$ and $1
\leq
m
\leq
n_2$, then
$$
L(a_1) \times L(a_2)=L(a_1+a_2).
$$
\endproclaim

The above proposition may be helpful in constructing  
complementary
series.

Let us now consider (U1). Certainly if $F=\Bbb C$, then 
there   is nothing to
prove since as we already mentioned, $D^u=(\Bbb 
C^\times)\hat{\ }$ and for 
$\chi \in D^u,\;\; u(\chi,n)=\chi \circ \det _n$.

For $F=\Bbb R$, $D^u \subseteq \GL (1,\Bbb R)\hat{\ } \cup 
\GL (2,\Bbb R)\hat{\
}$. Again, (U1) is evident if $\chi \in (\Bbb 
R^\times)\hat{\ }$. Speh
considered the remaining case of $u(\delta,n),\;\; \delta 
\in \GL (2,\Bbb
R)\hat{\ } \cap D^u$ [Sp2]. She proved unitarity using 
adelic methods.
Surprisingly,  it  seems that Gelfand and Graev were 
already aware of this
series of representations in the 1950s (see [Sp2, Remark 
1.2.2.] about [GfGr]). 

For nonarchimedean $F$, we have determined in [Td3] the 
representations 
$u(\delta,n)$ through the ideas presented  in \S 7. 
Unitarizability is proved there essentially along those 
lines. It is also
possible to prove unitarizability by the method of Speh as 
was done in
the appendix of [Td3]. Note that  here $D^u \cap
\GL (n,F)\hat{\ } \ne \varnothing $ for all $n \geq1.$

For $F=\Bbb C$, (U2) was proved by 
Stein in  [St]. In general, (U2) follows from (U0), using 
the   
irreducibility of representations $\nu^\alpha u(\delta,n)
\times
\nu^{-\alpha}u(\delta,n), 
0<\alpha<1/2$, obtained from   Proposition
9.1 and from the analytic properties of intertwining 
operators. There is
also another method in the nonarchimedean case presented 
in  [Bn2].

To prove (U3), one considers the $u(\delta,n)$'s as   
polynomials and 
proves the irreducibility of these polynomials. Here one 
uses the fact
that $R$ is a graded ring and that $u(\delta,n)$'s are 
homogeneous
elements. In the proof, one uses basic facts about the 
composition series
of generalized principal series representations (one does 
not need more
detailed information, such as that obtained from 
Kazhdan-Lusztig type
multiplicity formulas). It is a bit surprising that, 
although we do not
know how to write down the polynomials $u(\delta,n)$, we 
can nevertheless
carry out the proof. For proofs of (U3) see [Td3] and 
[Td2]. The
statement (U3) is obvious for $u(\delta,1)$ by
 Proposition 6.1.

Finally, let us return to (U0). Let  $P_n$ denote the
subgroup of $\GL (n,F)$ of all matrices with bottom row 
equal to $(0,\dotsc,0,1)$.
Already  Gelfand and Naimark noticed the importance of the 
statement 

\roster
\item"{(I)}" {If $\pi \in \GL (n,F)\hat{\ }$, then $\pi|P_n$
is irreducible $(n \in \Bbb N)$.}
\endroster

Actually, they proved the above statement  for $F=\Bbb C$, 
for the
representations that they expected to exhaust the unitary 
dual of $\GL (n,\Bbb
C)$.  Several people were aware that the above statement 
implies (U0). For a
written proof see  [Sh]. Proof of the implication is based 
on Mackey   theory and
Gelfand-Naimark models. 

Kirillov stated (I) as a theorem in  [Ki1] for F an  
archimedean
field. There he sketched a proof.  Vahutinskii's 
classification of
representations of $\GL (3,\Bbb R)\hat{\ }$ was based on 
the proof. Having
in mind a  correspondence obtained by Mackey theory 
$$\align
&\widehat{P}_n \rightarrow
\GL (n-1, F) \hat{\ } \cup \GL (n-2, F)\hat{\ } \cup 
\cdots \\
&\qquad \cup \GL (2, F) 
\hat{\ }  \cup \GL (1, F) \hat{\ } \cup \GL (0, F) \hat{\ },
\endalign
$$
the statement (I) would imply that one would have simpler  
realizations for 
representations from $\GL (n,F)\hat{\ 
}$. It is from this setup that the name Kirillov
model appears. In [Ki1] Kirillov's intention was to prove 
that $\pi|
P_n$ is operator irreducible (i.e., the commutator    
consists only of 
scalars), which is enough by Schur's lemma to see the 
irreducibility of
$\pi| P_n$. One takes any $T$ from the commutator of
$\pi|P_n$ and considers a distribution
$$
\Lambda_T:\vp \mapsto \text{Trace} \, (T\pi(\vp))
$$
on $\GL (n,F)$ which is  invariant for conjugations with 
elements from $P_n$, 
since $T$ is  $P_n$-intertwining.
 Now if 
$\Lambda_T$ is $\GL (n,F)$-invariant, then using the 
irreducibility of
$\pi$, it is not difficult to obtain that $T$ must be a 
scalar.  
Kirillov indicated that he proved in [Ki1]
that 
$\Lambda_T$ is $\GL (n,F)$-invariant (however, see below). 
This property of
the distribution is  especially easy to see when $\pi$
is a continuous finite-dimensional representation. Then the
distribution $\Lambda_T$ is  given by a continuous 
function which must be 
constant on $P_n$-conjugacy classes $(\Lambda_T$ is 
$P_n$-invariant). 
Since in $\GL (n,F)$, $\GL (n,F)$-conjugacy  classes 
contain dense
$P_n$-conjugacy classes,  $\Lambda_T$ must be $\GL 
(n,F)$-invariant.

Bernstein proved in  [Bn2] that for $F$ nonarchimedean, 
each $P_n$-invariant
distribution is $\GL (n,F)$-invariant. He proved (I) using 
essentially the
Kirillov's strategy.  Besides proving (U0) in the 
nonarchimedean case, he gave
a different proof of the implication (I) $\Rightarrow$ 
(U0). 

Bernstein states in  [Bn2] that Kirillov's proof of (I) 
for $F$ archimedean in
[Ki1] is incorrect, and he wrote that he himself had  an  
almost complete proof
(see [Bn2, p.\ 55]). In any case, there is no written 
complete proof of (I) in
the archimedean case now. We would say that  Kirillov's 
proof is incomplete
rather than incorrect. He failed to give a complete 
argument  that the
distribution  $\Lambda_T$ is $\GL (n,F)$-invariant. As was 
noted by Bernstein,
the tools used in Kirillov's paper do not seem to be 
sufficient for proving
(I). The distribution $\Lambda_T$  is a very special one. 
Actually (I) would
imply that it is a multiple of an irreducible character; 
so by Harish-Chandra's
regularity theorem,  it is locally $L^1$ and analytic on 
regular semisimple
elements. So if one proves that the (eigen) distribution  
$\Lambda_T$ is
locally $L^1$ and analytic on regular semisimple elements, 
one could apply
finite-dimensional argument. This may provide a strategy 
to prove (I) in the
archimedean case. This would be a longer proof, and there 
are also some
disadvantages in proving (U0) through (I). We shall say a 
few words about these
disadvantages. Before that, observe that there is an 
implicit proof of (U0) in
[Vo3]. 

We have mentioned that the approach to the unitary dual of 
$\GL (n,F)$
is expected to be applicable to $\GL (n)$ over a central 
division 
$F$-algebra $A$. As far as we understand, Vogan's 
description of the
unitary dual of $\GL (n,\Bbb H)$ confirms this. Here (U0) 
cannot be proved
through (I), simply because (I) is false in general in 
this case. The simplest
example can be obtained for $\GL (2,A)\;\;(A \ne F)$. It 
is not difficult to see
that there exist (irreducible) tempered representations 
which are
reducible when restricted to a nontrivial parabolic 
subgroup. Thus, it
may be more reasonable to search for proof of (U0) which 
works also for
division algebras. There are some candidates for it (see 
the last remark
at the end of this section).

After all, we have the following:

\proclaim{Theorem} Let $F$ be a locally compact 
nondiscrete field. Set
$$
B=\{u(\delta,n),\;\;  [\nu^\alpha u(\delta,n)] \times 
[\nu^{-\alpha}
u(\delta,n)],\;\;\delta \in D^u,\;\;n \in \Bbb N,\;\; 0 < 
\alpha <\tfrac 12 \}.
$$
Then
\roster
\item"{(i)}" if $\sigma_1, \ldots , \sigma_k \in B$,
then $
\sigma_1 \times \cdots \times \sigma_k \in \Irr ^u${\rm ;}
\item"{(ii)}" if
$\pi \in \Irr ^u$, then there exist 
$\sigma_1, \ldots , \sigma_m \in B$, unique up to a 
permutation, such 
that
$$
\pi=\sigma_1 \times \cdots \times \sigma_m.
$$ 
\endroster
\vskip-\smallskipamount
\vskip-1\belowdisplayshortskip
\nobreak
\endproclaim

We remind the reader once more that there is no written 
complete   proof yet of
(U0) in the archimedean case, but there is a complete 
proof  [Vo3] of the above
theorem in this case. 

\rem{Remarks} (1) We give in
[Td4] a concrete realization of the topological space $\GL 
(n,F)\hat{\
}$ when $F$ is nonarchimedean.

(2) The last theorem, together with  [Ka], implies that 
$\pi \in \SL (n,\Bbb
C)\hat{\ }$ is isolated if and only if $\pi$ is the 
trivial representation and
$n \ne 2$.

(3) Let $F$ be nonarchimedean. Let
$\rho \in \Irr ^u$ be a representation having a nontrivial 
  compactly
supported modulo center matrix coefficient. Then for $k \in 
\Bbb N$, the representation
$$
\nu^{(k-1)}\rho \times \nu^{(k-3)}\rho \times \cdots 
\times \nu^{-(k-1)}\rho
$$
has a unique square integrable subquotient which will be 
denoted  by
$\delta( \rho,k)$. In this way one obtains all $D^u$ 
(see [Ze] and [Jc1]). In
[Td4]
we have proved that
$\pi \in
\GL (n,F)\hat{\ }$
is isolated modulo center if and only if
$\pi$
equals some $u(\delta(
\rho
,k),m)$ with $k
\ne
2$ and $m
\ne
2$.

(4) One could try to prove (U0) for  $\GL (n)$ over local 
division algebras by proving first the following 
conjecture: Let $A$ be 
a central local simple algebra, let $S$ be the subgroup of 
the diagonal
matrices in $\GL (2,A)$, let $N$ be the subgroup of upper 
triangular
unipotent elements in $\GL (2,A)$, and let $\sigma$ be an 
irreducible 
unitary representation of $S$. Then 
$\text{Ind} ^{\GL (2,A)}_{SN}(\sigma)$
should be irreducible.  

(5) A proof of (U3) and (U4) for
$\GL (n)$ over a local  nonarchimedean division algebra is 
contained in 
[Td7]. 
\endrem


\Refs
\widestnumber\key{BnZe2$^{^{\text{\bf GL (n)}}}$}


{\tenpoint
The references that we  include here are directed more  to 
an inexperienced
reader in the field of the representation theory than to 
the experienced one.
This is the reason that we have classified them    into 
several   groups  (see
below).   We have tried to avoid too many very technical 
references, which are
very common in the field. We include a number of 
expository and survey papers
(they also omit many technicalities).  A more demanding 
reader can also  find a
choice of relevant references for further reading. We have 
omitted a number of
important references in order not to confuse a reader who 
is not very familiar
with the field.

Let us explain our classification of  the references  into
several groups. The classification is not very rigid. To
each reference we have attached a  superscript
that indicates the group where it belongs.  Here is a 
description of the groups.


\roster
\item"{{\bf sur}}"  
denotes the group of   survey and expository papers.


\item"{{\bf gen}}"
denotes the group of general references. 


\item"{{\bf GL (n)}}"
denotes  the group of  papers that are directly related
to the topic of our paper. They  are very useful for the
further understanding of the topics discussed in our
paper. The complete proofs that are omitted in this paper 
can be found in the papers from this group.


\item"{{\bf his}}"
denotes the group of the historically  important
references for the  topic of this paper.


\item"{{\bf Lan}}" 
denotes the  group of references directly related to the
Langlands program and groups $\GL (n)$. These references 
are very    often
related to the various correspondences predicted  by the 
Langlands
program for groups $\GL (n)$. Generally, these references
are  directed to more demanding and experienced readers.


\item"{{\bf oth}}" 
denotes the group of remaining useful references. They are
mostly related  to unitarity. 
\endroster
}

\ref\key A$^{^{\text{\bf sur}}}$
\by J. Arthur 
\paper Automorphic representations  and number theory
\inbook 1980 Seminar on Harmonic Analysis 
\pages 3--51
\publ Amer. Math. Soc. 
\publaddr Providence, RI
\yr 1981
\endref

\ref\key Bb$^{^{\text{\bf oth}}}$
\by D. Barbasch 
\paper The unitary dual for complex classical groups
\jour Invent. Math.
\vol 96
\yr 1989
\pages 103--176
\endref

\ref\key Bg$^{^{\text{\bf his}}}$
\by V. Bargmann 
\paper Irreducible unitary representations of
the Lorentz group
\jour Ann. of Math.
\vol 48
\yr 1947
\pages 568--640
\endref

\ref\key Bn1$^{^{\text{\bf oth}}}$
\by J. Bernstein 
\paper All
reductive $p$-adic groups are tame
\jour Funct. Anal. Appl.
\vol 8
\yr 1974
\pages 91--93
\endref

\ref\key Bn2$^{^{\text{\bf GL (n)}}}$
\bysame 
\paper P-invariant  distributions on GL\<$(N)$ and the
classification of unitary representations of GL\<$(N)$
(nonarchimedean case) \inbook Lie Group
Representations II, Proceedings, 
Univ. Maryland
1982--83  \pages 50--102
\publ Lecture Notes
in Math., vol. 1041, Springer-Verlag
\publaddr Berlin
\yr 1984
\endref

\ref\key BnZe1$^{^{\text{\bf GL (n)}}}$
\by J. Bernstein  and A. V. Zelevinsky 
\paper Representations of the group $GL (n,F)$,
where $F$ is a local nonarchimedean field
\jour Uspekhi Mat. Nauk.
\vol 31
\yr 1976
\pages 5--70
\endref

\ref\key BnZe2$^{^{\text{\bf GL (n)}}}$
\bysame 
\paper Induced representations of reductive $p$-adic 
groups. \RM I \jour Ann. Sci. \'{E}cole Norm. Sup. (4)
\vol 10
\yr 1977
\pages 441--472
\endref
 
\ref\key Bl$^{^{\text{\bf gen}}}$
\by A. Borel 
\book Linear algebraic groups 
\publ Benjamin
\publaddr New York
\yr 1969
\endref

\ref\key BlWh$^{^{\text{\bf gen}}}$
\by A. Borel and N. Wallach 
\book Continuous cohomology, discrete subgroups, and 
representations of  reductive groups 
\publ Princeton Univ. Press
\publaddr Princeton, NJ
\yr 1980
\endref

\ref\key Bu1$^{^{\text{\bf gen}}}$
\by N. Bourbaki 
\paper Groupes de Lie r\'eels compacts
\inbook Groupes et Alg\`ebres de Lie, chapter 9
\publ Masson
\publaddr Paris
\yr 1982
\endref

\ref\key Bu2$^{^{\text{\bf gen}}}$
\bysame 
\paper Mesure de Haar
\inbook Int\'egration
\bookinfo chapter 7
\publ Hermann
\publaddr Paris
\yr 1963
\endref

\ref\key Cy$^{^{\text{\bf oth}}}$
\by H. Carayol 
\paper  Repr\'esentations cuspidales du groupe lin\'eaire
\jour Ann. Sci. \'{E}cole Norm. Sup  (4) 
\vol 17
\yr 1984
\pages 191--226
\endref

\ref\key Ct$^{^{\text{\bf sur}}}$
\by P. Cartier 
\paper Representations  of  $p$-adic groups; a survey
\inbook Proc. Sympos. Pure Math., vol. 33 
\pages 111--155
\publ Amer. Math. Soc.
\publaddr Providence, RI 
\yr 1979
\endref

\ref\key Cs$^{^{\text{\bf gen}}}$
\by W. Casselman 
\paper Introduction to the theory of admissible 
representations of $p$-adic reductive groups 
\jour preprint
\endref

\ref\key CsMi$^{^{\text{\bf oth}}}$
\by W. Casselman and D. Mili\v{c}i\'{c} 
\paper  Asymptotic behavior of matrix coefficients  of
admissible representations
 \jour Duke Math. J.
\vol 49 
\yr 1982
\pages 869--930
\endref

\ref\key $\CL^{^{\text{\bf sur}}}$
\by L. Clozel 
\paper Progr\`{e}s r\'ecents vers la  classification  du
dual unitaire des groupes r\'eductifs r\'eels 
\jour
\paperinfo S\'eminaire Bourbaki, no. 681 (1987); 
Ast\'erisque {\bf 152--153}
\yr 1987 \pages 229--252
\endref

\ref\key DeKaVi$^{^{\text{\bf Lan}}}$
\by P. Deligne, D. Kazhdan, and M.-F. Vign\'eras 
\yr 1984
\paper Repr\'esentations des alg\`ebres centrales simples 
$p$-adiques 
\inbook Repr\'esentations des Groupes
R\'eductifs sur un Corps Local, by J.-N. Bernstein, 
P. Deligne, D. Kazhdan, and M.-F. Vign\'eras \publ
Hermann \publaddr Paris
\endref

\ref\key Di$^{^{\text{\bf gen}}}$
\by J. Dixmier 
\book Les $C^\ast$-algebras et leurs repr\'esentations
\publ Gauthiers-Villars
\publaddr Paris
\yr 1969
\endref

\ref\key Du1$^{^{\text{\bf sur}}}$
\by M. Duflo 
\paper Repr\'esentations de carr\'e int\'egrable des
groupes semi-simples r\'eels
\jour  S\'eminaire Bourbaki, no. 508 (1977-1978), Lecture
Notes in Math., vol. 710
\publ   Springer-Verlag
\publaddr Berlin, 1979
\endref

\ref\key Du2$^{^{\text{\bf oth}}}$
\bysame 
\paper Th\'eorie de Mackey pour les groupes alg\`ebriques
\jour Acta Math.
\vol 149
\yr 1982
\pages 153--213
\endref

\ref\key Fe$^{^{\text{\bf oth}}}$
\by J. M. G. Fell 
\paper Nonunitary dual space of groups
\jour Acta Math.
\vol 114
\yr 1965
\pages 267--310
\endref

\ref\key Gb1$^{^{\text{\bf sur}}}$
\by S. Gelbart 
\book Automorphic forms on adele groups
\publ Ann. of Math. Stud., vol. 83, Princeton
Univ. Press
\publaddr Princeton, NJ
\yr 1975
\endref

\ref\key Gb2$^{^{\text{\bf sur}}}$
\bysame 
\paper Elliptic curves and automorphic representations
\jour Adv. in Math.
\vol 21
\yr 1976
\pages 235--292
\endref 

\ref\key Gb3$^{^{\text{\bf sur}}}$
\bysame 
\paper An elementary introduction to the Langlands program
\jour Bull. Amer. Math. Soc.
\vol 10
\yr 1984
\pages 177--219
\endref 

\ref\key GfGr$^{^{\text{\bf his}}}$
\by I. M. Gelfand and M. I. Graev 
\paper Unitary representations of the real unimodular
group \jour Izv. Akad. Nauk
SSSR Ser. Mat.
\vol 17
\yr 1953
\pages 189--248
\endref

\ref\key GfGrPi$^{^{\text{\bf gen}}}$
\by I. M. Gelfand, M. Graev, and Piatetski-Shapiro
\book Representation theory and automorphic functions
\publ  Saunders
\publaddr Philadelphia, PA
\yr 1969
\endref

\ref\key GfKa$^{^{\text{\bf GL (n)}}}$
\by I. M. Gelfand and D. A. Kazhdan 
\yr 1974
\paper Representations of $GL (n,k)$
\inbook Lie Groups and Their Representations
\pages 95--118
\publ Halstead Press
\publaddr Budapest
\endref

\ref\key GfN1$^{^{\text{\bf his}}}$
\by I. M. Gelfand and M. A. Naimark 
\paper Unitary
representations of the Lorentz group
\jour Izv. Akad. Nauk SSSR Ser. Mat. 
\vol 11
\yr 1947
\pages 411--504
\afterall (Russian)
\endref

\ref\key GfN2$^{^{\text{\bf his}}}$
\bysame 
\book Unit\"are
Darstellungen der Klassischen Gruppen
{\rm (German translation of Russian publication from
1950)} 
\publ Akademie Verlag \publaddr Berlin
\yr 1957
\endref

\ref\key GfRa$^{^{\text{\bf his}}}$
\by I. M. Gelfand and D. A. Raikov 
\paper Irreducible unitary representations of locally
compact groups \jour Mat. Sb.
\vol 13 {\rm (55)}
\yr 1943
\pages 301--316
\endref

\ref\key Ha1$^{^{\text{\bf sur}}}$
\by Harish-Chandra
\paper  Harmonic analysis on semisimple Lie groups
\jour Bull. Amer. Math. Soc.
\yr 1970
\pages 529--551
\vol 76
\endref

\ref\key Ha2$^{^{\text{\bf sur}}}$
\bysame 
\paper Harmonic analysis on reductive $p$-adic groups
\inbook Proc. Sympos. Pure Math., vol. 26 
\pages 167--192
\publ Amer. Math. Soc.
\publaddr Providence, RI 
\yr 1973
\endref

\ref\key Ha3$^{^{\text{\bf gen}}}$
\bysame 
\book  Collected papers
\publ Springer-Verlag
\publaddr Berlin
\yr 1983
\endref

\ref\key He$^{^{\text{\bf Lan}}}$
\by G. Henniart 
\paper On the local Langlands conjecture for $GL (n)$\RM: 
the
cyclic case \jour Ann. of Math. (2)
\vol 123
\yr 1986
\pages 145--203
\endref

\ref\key Ho$^{^{\text{\bf oth}}}$
\by R. Howe 
\paper Tamely ramified supercuspidal
representations of $GL _n$
\jour Pacific J. Math.
\vol 73
\yr 1977
\pages 437--460
\endref

\ref\key HoMr$^{^{\text{\bf oth}}}$
\by R. Howe and C. C. Moore 
\paper Asymptotic properties of unitary representations
\jour J. Funct. Anal. 
\vol 32 
\yr 1979
\pages 72--96
\endref

\ref\key Jc1$^{^{\text{\bf GL (n)}}}$
\by H. Jacquet 
\yr 1977
\paper Generic representations
\inbook Non-Commutative Harmonic Analysis
\pages 91--101
\publ Lecture Notes in Math., vol. 587, Springer-Verlag
\publaddr Berlin
\endref

\ref\key Jc2$^{^{\text{\bf GL (n)}}}$
\bysame 
\paper On the residual spectrum of $GL (n)$
\inbook Lie Group
Representations II, Proceedings, Univ. of Maryland
1982--83  \pages 185--208
\publ Lecture Notes
in Math., vol. 1041, Springer-Verlag
\publaddr Berlin
\yr 1984
\endref

\ref\key JcL$^{^{\text{\bf Lan}}}$
\by H. Jacquet and R. P. Langlands 
\book Automorphic forms on $GL (2)$
\publ  Lecture Notes in Math., vol. 114, Springer-Verlag
\publaddr Berlin
\yr 1970
\endref

\ref\key Jn$^{^{\text{\bf oth}}}$
\by  C. Jantzen 
\paper Degenerate principal series for symplectic groups
\jour Mem. Amer.
Math. Society, no. 488
\publ Amer. Math. Soc.
\publaddr Providence, RI, 1993
\endref

\ref\key Ka$^{^{\text{\bf oth}}}$
\by D. A. Kazhdan 
\paper Connection of the
dual space of a group with the structure of its closed
subgroups \jour Funct.
Anal. Appl.
\vol 1
\yr 1967
\pages 63--65
\endref

\ref\key Ki1$^{^{\text{\bf GL (n)}}}$
\by A. A. Kirillov 
\paper Infinite dimensional representations of the
general linear group \jour Dokl.
Akad. Nauk SSSR
\vol 114
\yr 1962
\pages 37--39
\endref

\ref\key Ki2$^{^{\text{\bf gen}}}$
\bysame 
\book Elements of the theory of
representations
\publ Springer-Verlag
\publaddr New York
\yr 1976
\endref

\ref\key Kn$^{^{\text{\bf gen}}}$
\by A. W. Knapp 
\book Representation theory of semisimple groups
\publ Princeton Univ. Press
\publaddr Princeton, NJ
\yr 1986
\endref

\ref\key KnZu$^{^{\text{\bf sur}}}$
\by A. W. Knapp and G. J. Zuckerman 
\paper Classification theorems for representations of
semisimple Lie groups
\inbook Non-Commutative Harmonic Analysis
\pages 138--159
\publ Lecture Notes
in Math., vol.  587, Springer-Verlag
\publaddr Berlin
\yr 1977
\endref

\ref\key KuMy$^{^{\text{\bf Lan}}}$
\by P. Kutzko and A. Moy  
\paper On the local
Langlands conjecture in prime dimension
\jour Ann. of Math. (2)
\vol 121
\yr 1985
\pages 495--517
\endref

\ref\key L1$^{^{\text{\bf sur}}}$
\by R. P. Langlands 
\book Problems in the theory of automorphic forms
\publ  Lecture Notes in Math., vol. 170, Springer-Verlag
\pages 18--86
\publaddr Berlin
\yr 1970
\endref

\ref\key L2$^{^{\text{\bf gen}}}$
\bysame 
\paper On the
classification of irreducible representations of real
algebraic groups \inbook Representation Theory and
Harmonic Analysis on Semisimple Lie Groups 
\eds P. J. Sally,
Jr., and D. A. Vogan, Jr. 
\publ Amer. Math. Soc.
\publaddr Providence, RI 
\yr 1989
\endref

\ref\key Ma$^{^{\text{\bf his}}}$
\by F. Mautner 
\paper Spherical functions over $p$-adic fields. \RM I
\jour Amer. J. Math.
\vol 80
\yr 1958
\pages 441--457
\endref

\ref\key Mi$^{^{\text{\bf oth}}}$
\by D. Mili\v{c}i\'{c} 
\paper On $C^\ast$-algebras
with bounded trace
\jour Glasnik Mat.
\vol 8 
\yr 1973
\pages 7--21
\endref 

\ref\key M{\oe}Wd$^{^{\text{\bf oth}}}$
\by C. M{\oe}glin and J.-L. Waldspurger 
\paper Le spectre residuel de $GL (n)$
\jour Ann. Sci. \'{E}cole Norm. Sup. (4)
\vol 22
\yr 1989
\pages 605--674
\endref

\ref\key My$^{^{\text{\bf Lan}}}$
\by A. Moy 
\paper Local constants and the tame Langlands
correspondence \jour Amer. J. Math.
\vol 108
\yr 1986
\pages 863--930
\endref

\ref\key Ro$^{^{\text{\bf sur}}}$
\by F. Rodier 
\paper Repr\'esentations de $GL (n,k)$ o\`u $k$ est un
corps $p$-adique
\jour S\'eminaire Bourbaki, no. 587 (1982), Ast\'erisque
\vol 92--93
\yr 1982
\pages 201--218
\endref

\ref\key Sh$^{^{\text{\bf GL (n)}}}$
\by S. Sahi 
\paper On Kirillov's conjecture for archimedean fields
\jour Compositio Math.
\vol 72
\yr 1989
\pages 67--86
\endref

\ref\key SlTd$^{^{\text{\bf oth}}}$
\by P. J. Sally and M. Tadi\'{c} 
\paper Induced  representations and classifications for
$GSp(2, F)$ and $Sp(2,F)$
\jour M\'em. Soc. Math. France
\vol 52 \yr 1993
\endref

\ref\key Sd1$^{^{\text{\bf oth}}}$
\by F. Shahidi 
\paper Fourier transforms of intertwining operators and
Plancherel measure for $GL (n)$
 \jour Amer. J. Math.
\vol 106
\yr 1984
\pages 67--111
\endref

\ref\key Sd2$^{^{\text{\bf oth}}}$
\bysame 
\paper A proof of Langlands conjecture on Plancherel
measures; complementary series for $p$-adic groups
\jour Ann. of Math. (2)
\vol 132
\yr 1990
\pages 273--330
\endref

\ref\key Si1$^{^{\text{\bf oth}}}$
\by  A. Silberger 
\paper The Langlands quotient theorem for $p$-adic
groups
\jour Math. Ann.
\vol 236
\yr 1978
\pages 95--104
\endref

\ref\key Si2$^{^{\text{\bf gen}}}$
\bysame 
\book Introduction to harmonic analysis on reductive
$p$-adic groups \publ Princeton Univ. Press
\publaddr Princeton, NJ 
\yr 1979
\endref

\ref\key Sp1$^{^{\text{\bf oth}}}$
\by B. Speh 
\paper The unitary dual of
$GL (3,\Bbb R)$ and $GL (4,\Bbb R)$ 
\jour Math. Ann. 
\vol 258
\yr 1981
\pages 113--133
\endref

\ref\key Sp2$^{^{\text{\bf GL (n)}}}$
\bysame 
\paper Unitary
representations of $GL (n,\Bbb R)$ with nontrivial
$(\germ g
	,K)$-cohomology
\jour Invent. Math.
\vol 71
\yr 1983
\pages 443--465
\endref

\ref\key St$^{^{\text{\bf GL (n)}}}$
\by E. M. Stein 
\paper Analysis in matrix spaces and some new
representations of $\SL (N,\Bbb C)$
\jour Ann. of Math.
\vol 86
\yr 1967
\pages 461--490
\endref

\ref\key Td1$^{^{\text{\bf sur}}}$
\by M. Tadi\'{c} 
\paper  Unitary dual  of $p$-adic $GL (n)$, proof  of
Bernstein conjectures \jour Bull. Amer. Math. Soc.
\yr 1985
\vol 13
\pages 39--42
\endref

\ref\key Td2$^{^{\text{\bf GL (n)}}}$
\bysame 
\paper  Unitary  representations of general linear group
over  real and complex field
\jour preprint MPI/SFB 85-22 Bonn, 1985
\endref

\ref\key Td3$^{^{\text{\bf GL(n)}}}$
\bysame 
\paper  Classification of unitary representations in
irreducible representations of general  linear group
(nonarchimedean case) \jour Ann. Sci. \'{E}cole Norm.
Sup.  (4) \vol 19
\yr 1986
\pages 335--382
\endref

\ref\key Td4$^{^{\text{\bf oth}}}$
\bysame 
\paper Topology of unitary dual of nonarchimedean
$GL (n)$ \jour Duke Math. J.
\vol 55
\yr 1987
\pages 385--422
\endref

\ref\key Td5$^{^{\text{\bf oth}}}$
\bysame 
\paper On limits of characters of irreducible unitary 
representations \jour Glasnik Mat.
\vol 23 
\yr 1988
\pages 15--25
\endref

\ref\key Td6$^{^{\text{\bf oth}}}$
\bysame 
\paper Geometry of dual spaces of reductive  groups 
\RM(nonarchimedean case\/\RM) 
\jour J. Analyse Math.
\vol 51
\yr 1988
\pages 139--181
\endref

\ref\key Td7$^{^{\text{\bf oth}}}$
\bysame 
\paper Induced representations of $GL (n,A)$ for
$p$-adic division algebras $A$
\jour J. Reine Angew. Math.
\vol 405
\yr 1990
\pages 48--77
\endref

\ref\key Td8$^{^{\text{\bf oth}}}$
\bysame 
\paper On Jacquet modules of induced representations of
$p$-adic symplectic groups
\yr 1991
\inbook Harmonic Analysis on Reductive Groups
\pages 305--314
\publ Proceedings, Bowdoin College 1989,  Progr. 
Math., vol. 101, Birkh\"auser
\publaddr Boston, MA
\endref

\ref\key Vo1$^{^{\text{\bf gen}}}$
\by D. A. Vogan 
\book Representations of
real reductive groups 
\publ Birkh\"{a}user
\publaddr Boston, MA
\yr 1981
\endref

\ref\key Vo2$^{^{\text{\bf oth}}}$
\bysame
\paper Unitarizability of certain series of
representations \jour Ann. of Math. (2)
\vol 120
\yr 1984
\pages 141--187
\endref

\ref\key Vo3$^{^{\text{\bf GL(n)}}}$
\bysame 
\paper The unitary dual of $GL (n)$ over an archimedean
field
\jour Invent. Math.
\vol 82
\yr 1986
\pages 449--505
\endref

\ref\key Vo4$^{^{\text{\bf oth}}}$
\bysame 
\book Unitary
representations of reductive Lie groups 
\publ Princeton Univ. Press  
\publaddr Princeton, NJ
\yr 1987
\endref

\ref\key Wh$^{^{\text{\bf gen}}}$
\by N. Wallach 
\book Real reductive groups, \RM {vol. 1}
\publ Academic Press
\publaddr New York
\yr 1988
\endref

\ref\key Wr$^{^{\text{\bf gen}}}$
\by G.  Warner 
\book Harmonic analysis on semi-simple Lie groups. {\rm I, 
II}
\publ Springer-Verlag
\publaddr Berlin
\yr 1972
\endref

\ref\key We$^{^{\text{\bf gen}}}$
\by A. Weil 
\book Basic number theory
\publ Springer-Verlag
\publaddr New York
\yr 1974
\endref

\ref\key Ze$^{^{\text{\bf GL(n)}}}$
\by A. V. Zelevinsky 
\paper Induced representations of reductive
$p$-adic groups {\rm II}, On irreducible representations of
$GL (n)$ \jour Ann. Sci.
\'{E}cole Norm. Sup. (4)
\vol 13
\yr 1980
\pages 165--210
\endref
\endRefs

\enddocument